\documentclass[12pt,a4paper,dvips]{article}
\usepackage[cp1251]{inputenc}
\usepackage[english]{babel}
\usepackage{amsmath,amssymb,amsthm}
\usepackage[dvips,pdftex]{graphicx}
\usepackage{multirow}
\usepackage{rotating}
\usepackage{geometry,comment}
\usepackage[dvips,pdftex]{xcolor}
\usepackage{authblk}
\usepackage{epstopdf}

\newcommand{\R}{\mathbb{R}}
\newcommand{\Z}{\mathbb{Z}}

\newcommand{\SO}{\mathrm{SO}}
\newcommand{\SL}{\mathrm{SL}}
\newcommand{\PSL}{\mathrm{PSL}}
\newcommand{\SE}{\mathrm{SE}}
\newcommand{\SU}{\mathrm{SU}}

\newcommand{\sspan}{\mathrm{span}}

\newcommand{\Exp}{\mathrm{Exp}}
\newcommand{\ad}{\mathrm{ad}}
\newcommand{\Ad}{\mathrm{Ad}}
\newcommand{\const}{\mathrm{const}}
\newcommand{\VectField}{\mathrm{Vec}}
\newcommand{\sgn}{\mathrm{sgn}}
\newcommand{\id}{\mathrm{id}}
\newcommand{\diag}{\mathrm{diag}}
\newcommand{\Kil}{\mathrm{Kil}}
\newcommand{\type}{\mathrm{type}}
\newcommand{\M}{\mathcal{M}}
\newcommand{\p}{\bar{p}_3}
\newcommand{\st}{\sin{\tau}}
\newcommand{\ct}{\cos{\tau}}
\newcommand{\sht}{\sinh{\tau}}
\newcommand{\cht}{\cosh{\tau}}
\newcommand{\se}{\sin{(\tau \eta \p)}}
\newcommand{\ce}{\cos{(\tau \eta \p)}}
\newcommand{\argl}{\frac{t \eta p_3}{2 I_1}}
\newcommand{\stl}{\sin{\argl}}
\newcommand{\ctl}{\cos{\argl}}
\newcommand{\ImPart}{\mathrm{Im}}

\theoremstyle{definition}
\newtheorem{Def}{Definition}
\newtheorem{Remark}{Remark}

\theoremstyle{plain}
\newtheorem{Corollary}{Corollary}
\newtheorem{Lemma}{Lemma}
\newtheorem{Theorem}{Theorem}
\newtheorem{Proposition}{Proposition}

\newcommand*{\affaddr}[1]{#1} 
\newcommand*{\affmark}[1][*]{\textsuperscript{#1}}
\newcommand*{\email}[1]{\texttt{#1}}

\title{Symmetric Riemannian problem on the group of proper isometries of hyperbolic plane\footnote{This work is supported by the Russian Science Foundation
under grant 17-11-01387 and performed in Ailamazyan Program Systems
Institute of Russian Academy of Sciences.}}

\author{A.~V.~Podobryaev\affmark[1], Yu.~L.~Sachkov\affmark[2]\\
\affaddr{Program Systems Institute of RAS}\\
\email{\affmark[1]alex@alex.botik.ru, \affmark[2]yusachkov@gmail.com}
}

\begin{document}

\maketitle

\begin{abstract}
We consider the Lie group $\PSL_2(\R)$ (the group of orientation preserving isometries of the
hyperbolic plane) and a left-invariant Riemannian metric on this group with two equal eigenvalues
that correspond to space-like eigenvectors (with respect to the Killing form). For such metrics we
find a parametrization of geodesics, the conjugate time, the cut time and the cut locus. The
injectivity radius is computed. We show that the cut time and the cut locus in such Riemannian
problem converge to the cut time and the cut locus in the corresponding sub-Riemannian problem as
the third eigenvalue of the metric tends to infinity. Similar results are also
obtained for $\SL_2(\R)$.



\end{abstract}

\section{\label{section-introduction}Introduction}

The Riemannian problem is a problem of finding shortest arcs connecting two arbitrary points of
a Riemannian manifold. If this manifold is a homogeneous space of a group $G$ and Riemannian metric
is invariant under the action of $G$, then we can consider only geodesics starting at the fixed
point.
So, the description of the shortest arcs is equivalent to the description of the cut time
and cut points of such geodesics. Recall that the cut time is the time of loss of optimality
of a geodesic. The cut point is a geodesic's point that corresponds to the cut time. The cut locus is
the union of cut points of all geodesics starting at the fixed point.

The main result of this paper is the description of the cut locus and the cut time for the series of
symmetric Riemannian problems on the group $\PSL_2(\R)$.

In our previous work~\cite{podobryaev-sachkov-dan} we show that this series joins with the series of
the symmetric Riemannian problems on $\SO_3$ extended by two sub-Riemannian problems on
$\PSL_2(\R)$ and $\SO_3$. Now we give a brief description of known results on this extended series of
Riemannian and sub-Riemannian problems.

Let $X$ be a two-dimensional Riemannian manifold of a constant non zero curvature
$\kappa = -1 \ \text{or} \ 1$. So, $X$ is the hyperbolic (Lobachevsky) plane $\Lambda^2$ or the sphere
$S^2$.
Let $G$ be the group of isometries preserving the orientation of $X$, i.e., $G$ is $\PSL_2(\R)$
or $\SO_3$ respectively.

Consider $G$ as the bundle $\mathrm{S}X$ of unit tangent vectors to $X$.
This bundle is a weakly symmetric space $(G \times \SO_2) / \SO_2$. The second multiplier $\SO_2$ acts
on the bundle of unit tangent vectors to $X$ by rotations by the same angle in all tangent spaces. The
stabilizer $\SO_2$ is embedded into the direct product in the anti-diagonal way. Weakly symmetric
spaces were introduced by A.~Selberg~\cite{selberg} and
$(\PSL_2(\R) \times \SO_2) / \SO_2$ is Selberg's first original example. We consider a $(G \times
\SO_2)$-invariant Riemannian metric on $G$. In other words it is a left-invariant Riemannian metric on
$G$ which is a lift of a Riemannian metric on $X$. Such a metric is determined by three eigenvalues
$I_1 = I_2, I_3$ of the restriction of the metric to the tangent space at the identity.
We call the corresponding Riemannian problem the symmetric Riemannian problem.

Now fix the distribution of two-dimensional planes in $G$ that are orthogonal to fibres of the
projection from $\mathrm{S}X$ onto $X$. Consider the sub-Riemannian metric on $G$
defined by this distribution and the restriction of the Killing form to this distribution.
The sub-Riemannian problem is a problem of finding shortest arcs of sub-Riemannian geodesics.
The sub-Riemannian problems for $\PSL_2(\R)$ and $\SO_3$ were considered by V.~N.~Berestovskii and
I.~A.~Zubareva~\cite{berestovskii-zubareva-sl2, berestovskii-zubareva-so3, berestovskii},
and by U.~Boscain and F.~Rossi~\cite{boscain-rossi}.

In this paper we consider a series of Riemannian problems on the group of isometries of
the hyperbolic plane with $I_1 = I_2$. The cut locus and the equations for the cut time are found. It
turns out that the Riemannian problem approximates the sub-Riemannian one as $I_3 \rightarrow \infty$.
This means that the parametrization of geodesics, the conjugate time, the conjugate locus, the cut
time, the cut locus of the Riemannian problem converge to the same objects in the sub-Riemannian one.
We have achieved
similar results for $\SO_3$ (the group of isometries of a sphere) in~\cite{podobryaev-sachkov}.

Table~\ref{table-results} presents a summary of known results.
Here we use the following notation:
$$
\eta = \kappa\frac{I_1}{I_3} - 1,
$$
$Z$ is the set of all central symmetries of $X$, $R_{\eta}$ is the interval of some rotations of $X$
around the fixed point. This interval depends on the parameter $\eta$ and converges to the circle $R$
of all rotations of $X$ around the fixed point as $\eta \rightarrow -1$
(equivalent to $I_3 \rightarrow \infty$).

\begin{table}[h]
\label{table-results}
\caption{Summary of known results. For $\eta < -1$ we have Riemannian problem on $\PSL_2(\R)$,
for $\eta > -1$ we have Riemannian problem on $\SO_3$. Limits $\eta \rightarrow -1 \pm 0$
correspond to sub-Riemannian problems on $\PSL_2(\R)$ and $\SO_3$ respectively.}
\medskip
\begin{tabular}{|c|c|l|l|}
\hline
$\eta$ & $I_1, I_3$ & Closure & Reference \\
& & of the cut locus & \\
\hline
$\eta \leqslant -\frac{3}{2}$ & $2I_1 \geqslant I_3$ & $Z$ &
\multirow{2}{*}{a result of this paper} \\
\cline{1-3}
$-\frac{3}{2} < \eta < -1$ & $2I_1 < I_3$ & $Z \cup R_{\eta}$ & \\
\hline
$\eta \rightarrow -1-0$ & $I_3 \rightarrow \infty$ & $Z \cup R$ &
V.~N.~Berestovskii~\cite{berestovskii} \\
\hline
$\eta \rightarrow -1+0$ & $I_3 \rightarrow \infty$ & $Z \cup R$ &
V.~N.~Berestovskii, I.~A.~Zubareva~\cite{berestovskii-zubareva-so3} \\ & & & U.~Boscain,
F.~Rossi~\cite{boscain-rossi} \\
\hline
$-1 < \eta < -\frac{1}{2}$ & $2I_1 < I_3$ & $Z \cup R_{\eta}$ & \multirow{2}{*}{A.~V.~Podobryaev,
Yu.~L.~Sachkov~\cite{podobryaev-sachkov}} \\
\cline{1-3}
$-\frac{1}{2} \leqslant \eta$ & $2I_1 \geqslant I_3$ & $Z$ & \\
\hline
\end{tabular}
\end{table}

Also one can consider the Euclidian plane $X = E^2$, a two-dimensional manifold of constant zero curvature
and the group $G = \SE_2$ of isometries preserving the orientation of $E^2$. The answer in the
corresponding Riemannian problem is unknown. The sub-Riemannian problem on the upper defined
distribution is not completely controllable. But there is a result of the sub-Riemannian problem
for another distribution of two-dimensional planes tangent to the fibres of the projection from
$\mathrm{S}X$ onto $X$. This sub-Riemannian problem models a vehicle on $X = E^2$ that can go forward
and can rotate. This problem was solved by the second co-author in the series of
papers~\cite{sachkov-moiseev-se-1} (in collaboration with I.~Moiseev), \cite{sachkov-se-2},
\cite{sachkov-se-3}. The cut locus is the union of the set of all central symmetries and a part of
M\"{o}bius strip with unknown geometric sense.

Now introduce our plan of investigation of the cut locus on $\PSL_2(\R)$:
\begin{enumerate}
\item parametrization of geodesics via the Pontryagin maximum principle~\cite{pontryagin, agrachev-sachkov};
\item description of the group of symmetries of the exponential map;
\item description of the Maxwell strata and the Maxwell time that correspond to the symmetry group
    of the exponential map;
\item finding the first conjugate time;
\item it turns out that the first conjugate time is greater than (or equal to) the Maxwell time
(corresponding to the symmetries), and the exponential map is a diffeomorphism of the set bounded by
the first Maxwell time in the pre-image of the exponential map to an open dense subset of $G$. That
is why the first Maxwell time turns out to be the cut time. Then we describe the global structure of
the cut locus.
\end{enumerate}

This scheme of investigation of the global optimality of extremals first appears in works
of the second co-author on the generalized Dido problem~\cite{sachkov-didona1, sachkov-didona2}.

The structure of this paper corresponds to the above items
(Sections~\ref{section-geodesics-parametrization}--\ref{section-cut-set}).

In Section~\ref{section-injection-radius} the injectivity radius of the considered metric is computed
(depending on the parameters $I_1$ and $I_3$).

Section~\ref{section-sl2} contains results on the similar Riemannian problem on the group $\SL_2(\R)$
(an answer in the corresponding sub-Riemannian problem was achieved by U.~Boscain and
F.~Rossi~\cite{boscain-rossi}, a complete proof was given by V.~N.~Berestovskii and
I.~A.~Zubareva~\cite{berestovskii-zubareva-sl2}, besides such sub-Riemannian problem was considered by
E.~Grong and A.~Vassil'ev~\cite{grong-vasil'ev}). Note that we got a result in the similar
Riemannian problem on $\SU_2$ \cite{podobryaev-sachkov}, while the corresponding sub-Riemannian
problem
on $\SU_2$ was considered by D.-Ch.~Chang, I.~Markina and
A.~Vassil'ev~\cite{chang-markina-vasil'ev}.

Section~\ref{section-sub-riemannian} deals with the Riemannian approximation of the sub-Riemannian
problem as $I_3 \rightarrow \infty$.

\section{\label{section-geodesics-parametrization}Parametrization of geodesics}

\subsection{\label{section-def}Definitions and notation}
Let $G = \PSL_2(\R)$ and let $\mathfrak{g}$ be the corresponding Lie algebra. Consider a basis of the Lie algebra  $e_1, e_2, e_3 \in \mathfrak{g}$ such that the Killing form and the Riemannian metric have the matrices $\diag(1, 1, -1)$ and $\diag(I_1, I_2, I_3)$ respectively. Next we consider the case of  $I_1 = I_2, I_3 > 0$. By
$$
\eta = -\frac{I_1}{I_3} - 1 < -1
$$
denote a parameter of the Riemannian metric. This parameter measures prolateness of small spheres.
We identify $\mathfrak{g}$ with $\mathfrak{g}^*$ via the Killing form.
Assume that this identification takes $e_1, e_2, e_3 \in \mathfrak{g}$ to a basis $\varepsilon_1, \varepsilon_2, \varepsilon_3 \in \mathfrak{g}^*$.
Let $p = p_1 \varepsilon_1 + p_2 \varepsilon_2 + p_3 \varepsilon_3 \in \mathfrak{g}^*$.

Introduce the following notation:
$$
\Kil(p) = p_1^2 + p_2^2 - p_3^2, \qquad |p| = \sqrt{|\Kil(p)|}, \qquad \type(p) = \sgn(-\Kil(p)),
$$
where $\Kil(p)$ is the value of quadratic Killing form on a covector $p$. Recall that $p$ is called
time-like, light-like or space-like if $\type(p)$ is equal to $1$, $0$ or $-1$ respectively.

Assume that all geodesics have an arclength parametrization by a parameter $t$ (called a time).
For $|p| \neq 0$ define
$$
\bar{p} = \frac{p}{|p|}, \qquad \tau(p) = \frac{t |p|}{2 I_1}.
$$

By $R_{v, \varphi}$ denote the rotation of a three-dimensional oriented Euclidean space around the axis $\sspan\{v\}$ by the angle $\varphi$ in the positive direction.

\subsection{\label{section-optimal-control}Optimal control problem}

Consider the problem of finding shortest arcs of the Riemannian metric as an optimal control problem~\cite{agrachev-sachkov}:
\begin{equation}
\label{eq-optimal-control-problem}
\begin{aligned}
\begin{array}{rcl}
\dot{Q} = Q \Omega, & \quad & \Omega = u_1 e_1 + u_2 e_2 + u_3 e_3 \in \mathfrak{g}, \\
Q \in G, & \quad & (u_1, u_2, u_3) \in \R^3, \\
Q(0) = \id, & \quad & Q(t_1) = Q_1, \\
\end{array} \\
\frac{1}{2} \int_0^{t_1}{(I_1 u_1^2 + I_2 u_2^2 + I_3 u_3^2) \ dt} \rightarrow \min,
\end{aligned}
\end{equation}
where $u = (u_1, u_2, u_3)$ is a control.
Minimization of the Riemannian length is equivalent to minimization of this energy functional
due to the Cauchy-Schwartz inequality (with a fixed terminal time $t_1$).

\subsection{\label{section-geodesics-eq}Equations of geodesics}

The following theorem gives a parametrization of geodesics.

{\Theorem
\label{th-geodesics-parametrization}
A geodesic $Q(t)$ starting at the identity and having an initial momentum
$p = p_1 \varepsilon_1 + p_2 \varepsilon_2 + p_3 \varepsilon_3$ $(\text{where} \quad
\frac{p_1^2}{I_1} + \frac{p_2^2}{I_2} + \frac{p_3^2}{I_3} = 1)$
is a product of two one-parameter subgroups:
\begin{equation}
\label{eq-geodesic-parametrization}
Q(t) = \exp \left(\frac{tp}{I_1}\right) \exp \left(\frac{t \eta p_3 e_3}{I_1} \right).
\end{equation}
}
\medskip

\emph{Proof.}
Geodesics are extremals of the optimal control problem~(\ref{eq-optimal-control-problem}). Apply the Pontryagin maximum principle~\cite{pontryagin}. Consider the trivialization of the cotangent bundle
$\tau: G \times \mathfrak{g}^* \rightarrow T^* G$ via the $G$-action: $\tau(g, \alpha) = dL_g^* \alpha$, where $L_g : G \rightarrow G$ is the left shift by $g \in G$, and $\alpha \in \mathfrak{g}^*$.

The Hamiltonian of the Pontryagin maximum principle reads as
$$
H_u^\nu(p) = u_1 p_1 + u_2 p_2 - u_3 p_3 + \frac{\nu}{2}(I_1 u_1^2 + I_2 u_2^2 + I_3 u_3^2),
$$
where $\nu \leqslant 0$.
For an extremal control $\tilde{u}(t)$ for almost any time
$H_{\tilde{u}(t)}^\nu(p(t)) = \max_u H_u^\nu (p(t))$. As usual in Riemannian problems $\nu = 0$ implies $p = 0$ in contradiction with the condition of Pontryagin maximum principle
of non-triviality of the pair $(\nu, p)$. This pair is defined up to a positive multiplier. So we can set $\nu = -1$. Then
$$
\tilde{u}_1(t) = \frac{p_1(t)}{I_1}, \qquad \tilde{u}_2(t) = \frac{p_2(t)}{I_2}, \qquad \tilde{u}_3(t) = -\frac{p_3(t)}{I_3}.
$$
The maximized Hamiltonian is
$$
H(p) = H^{-1}_{\tilde{u}(t)}(p) = \frac{1}{2} \left(\frac{p_1^2}{I_1} + \frac{p_2^2}{I_2} + \frac{p_3^2}{I_3}\right).
$$
The corresponding Hamiltonian system reads as
\begin{equation*}
\left\{
\begin{aligned}
\begin{array}{ccl}
\dot{Q}(t) & = & Q(t) \ \Omega(t),\\
\dot{p}(t) & = & (\ad^* \ \Omega(t))p(t), \\
\end{array}
\end{aligned}
\right.
\end{equation*}
where $\Omega(t) = \tilde{u}_1(t) e_1 + \tilde{u}_2(t) e_2 + \tilde{u}_3(t) e_3 \in \mathfrak{g}$.
We call the first equation \emph{the horizontal part} and the second one \emph{the vertical part} of the Hamiltonian system.

It is easy to see that in the coordinates $p_1, p_2, p_3$ the equations of the vertical part are
as follows:
\begin{equation*}
\left\{
\begin{aligned}
\begin{array}{ccl}
\dot{p}_1(t) & = & -p_2(t) p_3(t) \frac{I_2 + I_3}{I_2 I_3}, \\
\dot{p}_2(t) & = & p_1(t) p_3(t) \frac{I_1 + I_3}{I_1 I_3}, \\
\dot{p}_3(t) & = & p_1(t) p_2(t) \frac{I_1 - I_2}{I_1 I_2}. \\
\end{array}
\end{aligned}
\right.
\end{equation*}

When $I_1 = I_2$ the solution is
\begin{equation}
\label{eq-geodesic-vertical}
p(t) = R_{e_3, -\frac{t \eta p_3}{I_1}} p, \qquad p(0) = p.
\end{equation}
In invariant notation
$$
p(t) = \Ad \exp \left( -\frac{t \eta p_3 e_3}{I_1} \right) p.
$$
Note that
$$
\Omega(t) = \frac{p_1(t)}{I_1} e_1 + \frac{p_2(t)}{I_2} e_2 - \frac{p_3(t)}{I_3} e_3 =
\frac{1}{I_1} (p(t) + \eta p_3 e_3) =
\frac{1}{I_1} \Ad \exp \left( -\frac{t \eta p_3 e_3}{I_1} \right) (p + \eta p_3 e_3).
$$
Compute derivative of $Q(t)$, see~(\ref{eq-geodesic-parametrization}):
\begin{equation*}
\begin{aligned}
\begin{array}{ccl}
\dot{Q}(t) & = & dL_{\exp{\frac{tp}{I_1}}} dR_{\exp{\frac{t \eta p_3 e_3}{I_1}}} \left( \frac{p}{I_1} \right) +
dL_{\exp{\frac{tp}{I_1}}} dL_{\exp{\frac{t \eta p_3 e_3}{I_1}}} \left( \frac{\eta p_3 e_3}{I_1} \right) \\
 & = & dL_{\exp{\frac{tp}{I_1}} \exp{\frac{t \eta p_3 e_3}{I_1}}}
\left(
\frac{1}{I_1} dL_{\exp{(-\frac{t \eta p_3 e_3}{I_1})}} dR_{\exp{(\frac{t \eta p_3 e_3}{I_1})}}
(p + \eta p_3 e_3)
\right),\\
\end{array}
\end{aligned}
\end{equation*}
where $L_g$ and $R_g$ are the left and right shifts by $g \in G$ respectively. From the formula for $\Omega(t)$ it follows that the last expression is equal to $Q(t) \Omega(t)$. So, $Q(t)$ satisfies the horizontal part of the Hamiltonian system.
\qquad$\Box$
\medskip

{\Remark
The solution of the vertical part of the Hamiltonian system takes the form
$\bar{p}(\tau) = R_{e_3, -2 \tau \eta \p} \bar{p}$, if $|p| \neq 0$.
}
\medskip

{\Remark
The Killing form is a Casimir function on $\mathfrak{g}^*$. Thus, $\type(p) \equiv \const$, i.e.,
type of a covector is an integral of the Hamiltonian system.
}
\medskip

\subsection{\label{section-model}Model of $\PSL_2(\R)$}

Here we describe a model of $\PSL_2(\R)$ in which we produce computations and draw figures.
Consider the group $\SU_{1, 1}$, which can be realized as the group of unit norm split-quaternions
$$
\SU_{1, 1} = \{ q_0 + q_1 i + q_2 j + q_3 k \ | \ q_0^2 - q_1^2 - q_2^2 + q_3^2 = 1, \ q_0, q_1, q_2, q_3 \in \R\}.
$$
The product rule of split-quaternions is distributive and satisfies the following conditions:
\begin{equation*}
i^2 = j^2 = 1, \qquad k^2 = -1, \qquad
ij = -k, \qquad
jk = i, \qquad
ki = j.
\end{equation*}
It is well known (see for example \cite{vinberg-onishchik}) that there is an isomorphism
$$
\psi: \SL_2(\R) \rightarrow \SU_{1, 1}, \quad
\psi
\left(
\begin{array}{ll}
a & b \\
c & d
\end{array}
\right)
= \frac{a+d}{2} + \frac{a-d}{2} i + \frac{b+c}{2} j + \frac{c-b}{2} k, $$
where $ad-bc = 1$, $a, b, c, d \in \R$.

Consider a projection of $\SU_{1, 1}$ onto a three-dimensional real space with coordinates $q_1, q_2, q_3$. The condition
$$
q_3^2 - q_1^2 - q_2^2 = 1 - q_0^2 \leqslant 1
$$
implies that the image of $\SU_{1, 1}$ is a domain between two cups of the hyperboloid defined by the equation
$q_3^2 - q_1^2 - q_2^2 = 1$. For fixed $q_1, q_2, q_3$ (such that $q_3^2 - q_1^2 - q_2^2 \neq 1$)
there are two possibilities for the value of $q_0$. Hence, the group $\SU_{1, 1}$ is the union of two such domains with identified boundary points (which correspond to $q_0 = 0$). The group $\SU_{1, 1}$ is homeomorphic to an open solid torus.

The group $\PSL_2(\R) \backsimeq \SU_{1, 1}/\{\pm 1\}$ can be seen as the domain between the cups of the  hyperboloid with identified opposite points on the cups of the hyperboloid: $(q_1, q_2, q_3) \sim (-q_1, -q_2, -q_3)$.

\begin{figure}[h]
\caption{Model of $\SU_{1, 1}$.}
\label{pic-su11model}
\medskip
\center{\includegraphics[width=0.5\linewidth]{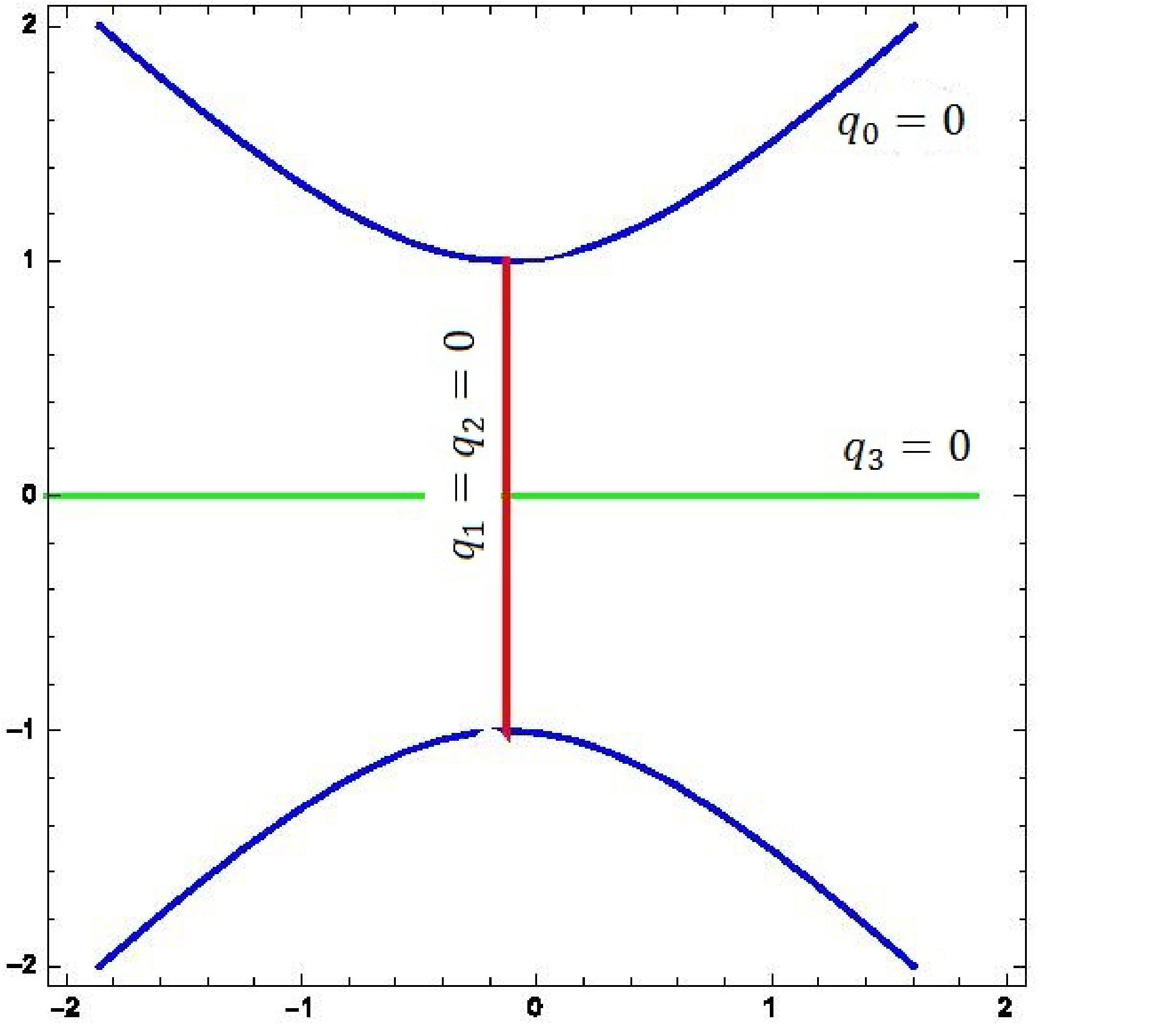}}
\end{figure}

As we will see below, the Maxwell strata and the cut locus are invariant under rotations around axis $q_3$ because of the symmetry of the Riemannian metric ($I_1 = I_2$).
So, we will draw all required sets on the plane $q_1, q_3$ (see Figure~\ref{pic-su11model}).

The group $\SU_{1,1} / \{ \pm E \}$ is the group of M\"{o}bius transformations of complex numbers that preserve the unit disk. In other words, it is the group of orientation preserving
isometries of the Poincar\'{e} disk model of the hyperbolic geometry. The split-quaternion $q_0 + q_1 i + q_2 j + q_3 k$ corresponds to the M\"{o}bius transformation
\begin{equation}
\label{eq-drlin}
z \mapsto \frac{(q_0 + q_3 i) z + (q_1 + q_2 i)}{(q_1 - q_2 i) z + (q_0 - q_3 i)}, \qquad z \in \mathbb{C}.
\end{equation}

It is well known (see Appendix) that an orientation preserving isometry of the hyperbolic plane is a product of two reflections with respect to lines.
There are three types of orientation preserving isometries: elliptic, parabolic and hyperbolic ones.
These types correspond to pairs of lines that are intersecting, parallel one to another
(the intersection point is on the absolute) or ultra parallel one to another (not intersecting)
respectively.

The following proposition provides a geometric interpretation of some subsets of $\SU_{1, 1}$.

{\Proposition
Consider the projection
$\Pi : \SU_{1, 1} \rightarrow \SU_{1, 1} / \{ \pm E \} \backsimeq \PSL_2(\R)$, then: \\
$(1)$ $\Pi \{q \in \SU_{1, 1} \ | \ q_1 = q_2 = 0 \}$ is the set of rotations around the center of the Poincar\'{e} disk model; \\
$(2)$ $\Pi \{q \in \SU_{1, 1} \ | \ q_0 = 0 \}$ is the set of central symmetries (reflections in points);\\
$(3)$ $\Pi \{q \in \SU_{1, 1} \ | \ q_3 = 0 \}$ is the set of hyperbolic isometries that is defined by a sheaf of ultra parallel lines that is symmetric with respect to a diameter of
the Poincar\'{e} disk model.
}
\medskip

\emph{Proof.}
(1) Obviously, the corresponding M\"{o}bius transformation is the multiplication by $(q_0 + q_3 i)^2$. That is rotation around zero by the angle $2 \arg{(q_0 + q_3 i)}$. \\
(2) If $z$ is a fixed point of the M\"{o}bius transformation~(\ref{eq-drlin}), then
\begin{equation}
\label{eq-fixpoint}
(q_1 - q_2i) z^2 - 2 q_3 i z - (q_1 + q_2 i) = 0.
\end{equation}
One of the two solutions $\frac{i (q_3 \pm 1)}{q_1 - q_2i}$ is inside of the unit disk. Hence, we have an elliptic isometry (rotation). Compute derivative of transformation~(\ref{eq-drlin}) at the fixed point:
$$
\frac{q_3^2 - (q_1^2 + q_2^2)}{(q_1 - q_2 i)^2 z^2 - 2 q_3 i (q_1 - q_2 i) z - q_3^2} =
\frac{1}{(q_1 - q_2 i) ((q_1 - q_2 i) z^2 - 2 q_3 i z) - q_3^2}.
$$
But $z$ satisfies equation~(\ref{eq-fixpoint}), then
$$
\frac{1}{(q_1 - q_2i)(q_1 + q_2i) - q_3^2} = -1.
$$
This implies that the transformation is the reflection in the point $z$ (the central symmetry). \\
(3) Fixed points of the M\"{o}bius transformation can be found from the equation
$$
(q_1 - q_2i) z^2 - (q_1 + q_2i) = 0.
$$
There are two opposite solutions with the same absolute value that is equal to $1$.
Thus, the transformation is a hyperbolic isometry. The corresponding sheaf of ultra parallel lines is symmetric with respect to the diameter connecting
the two fixed points.
\qquad$\Box$
\medskip

\subsection{\label{section-exponential-map}Exponential map}

{\Corollary
\label{crl-geodesics-parametrization-in-model}
A geodesic starting at the identity of $\SU_{1, 1}$ with an initial momentum $p$
has the following arclength parametrization: \\
$(1)$ for a time-like covector $p \quad (p_3^2 - p_1^2 - p_2^2 > 0)$
\begin{equation}
\label{eq-time-geodesic}
\begin{aligned}
\begin{array}{ccl}
q_0^e(\tau) & = & \ct \ce - \p \st \se, \\
\left(
  \begin{array}{l}
     q_1^e(\tau )\\
     q_2^e(\tau )
  \end{array}
\right) & = &  \st R_{e_3, -\tau \eta \p}
\left(
  \begin{array}{l}
     \bar{p}_1\\
     \bar{p}_2
  \end{array}
\right), \\
q_3^e(\tau ) & = & \ct \se + \p \st \ce,
\end{array}
\end{aligned}
\end{equation}
$(2)$ for a light-like covector $p \quad (p_3^2 - p_1^2 - p_2^2 = 0)$
\begin{equation}
\label{eq-light-geodesic}
\begin{aligned}
\begin{array}{ccl}
q_0^p(t) & = & \ctl - \frac{t}{2 I_1} p_3 \stl, \\
\left(
  \begin{array}{l}
     q_1^p(t)\\
     q_2^p(t)
  \end{array}
\right) & = &  \frac{t}{2 I_1} R_{e_3, -\argl}
\left(
  \begin{array}{l}
     p_1\\
     p_2
  \end{array}
\right), \\
q_3^p(t) & = & \stl + \frac{t}{2 I_1}p_3 \ctl,
\end{array}
\end{aligned}
\end{equation}
$(3)$ for a space-like covector $p \quad (p_3^2 - p_1^2 - p_2^2 < 0)$
\begin{equation}
\label{eq-space-geodesic}
\begin{aligned}
\begin{array}{ccl}
q_0^h(\tau) & = & \cht \ce - \p \sht \se, \\
\left(
  \begin{array}{l}
     q_1^h(\tau )\\
     q_2^h(\tau )
  \end{array}
\right) & = & \sht R_{e_3, -\tau \eta \p}
\left(
  \begin{array}{l}
     \bar{p}_1\\
     \bar{p}_2
  \end{array}
\right), \\
q_3^h(\tau ) & = & \cht \se + \p \sht \ce.
\end{array}
\end{aligned}
\end{equation}
}
\medskip

\emph{Proof.} Let $p = p_1 \frac{i}{2} + p_2 \frac{j}{2} + p_3 \frac{k}{2} \in \mathfrak{su}_{1, 1}$ be the orthonormal (with respect to the Killing form) decomposition of the vector $p$. Then
$$
p^n = \frac{\Kil(p)^{[\frac{n}{2}]}}{2^n} p^{(n \mod 2)}.
$$
Consider the exponential map from the Lie algebra to the Lie group
$ \exp: \mathfrak{su}_{1, 1} \rightarrow \SU_{1, 1}.$ It follows that
$$
\exp(p) =
\left\{
\begin{array}{rcc}
\cos{(\frac{|p|}{2})} + \sin{(\frac{|p|}{2})} (\bar{p}_1 i + \bar{p}_2 j + \bar{p}_3 k), & \text{for} & \Kil(p) < 0, \\
1 + \frac{1}{2}(p_1 i + p_2 j + p_3 k), & \text{for} & \Kil(p) = 0, \\
\cosh{(\frac{|p|}{2})} + \sinh{(\frac{|p|}{2})} (\bar{p}_1 i + \bar{p}_2 j + \bar{p}_3 k), & \text{for} & \Kil(p) > 0. \\
\end{array}
\right.
$$
It remains to calculate the product of the expressions of the two one-parametric subgroups from Theorem~\ref{th-geodesics-parametrization}.
\qquad$\Box$
\medskip

We will skip the upper index of the functions $q_0, q_1, q_2, q_3$ when we formulate a general statement for them.

{\Remark
The image of a geodesic under the projection $\Pi : \SU_{1,1} \rightarrow \PSL_2(\R)$ is a geodesic. Inversely any geodesic in $\PSL_2(\R)$ lifts to a geodesic in $\SU_{1, 1}$.
}
\medskip

{\Def
Let
$C=\{p \in \mathfrak{g}^* \ | \ H(p) = 1/2 \}$
be the level surface of the Hamiltonian (an ellipsoid).
Initial momenta from $C$ correspond to unit initial velocities of geodesics
(i.e., the arclength parametrization of geodesics).
}

{\Def
\label{def-exponential-map}
\emph{The exponential map} (for the Riemannian problem) is the map
$$
\Exp : C \times \R_{+} \rightarrow G, \qquad
\Exp(p, t) = \pi \circ e^{t\vec{H}}(\id, p),
$$
where $p \in C$, $t \in \R_+$ and $e^{t\vec{H}}$ is the flow of the Hamiltonian vector field
$\vec{H} \in \VectField(T^*G)$, and $\pi : T^*G \rightarrow G$ is the projection
of the cotangent bundle to the base.
}
\medskip

The exponential map defines the arclength parametrization of geodesics.

{\Remark
\label{rem-exponential-map-is-real-analitic}
The exponential map is real analytic, since the Hamiltonian $H$ and the Hamiltonian vector field $\vec{H}$ are real analytic.
}
\medskip

\section{\label{section-exponential-map-symmetries}Symmetries of exponential map}

In this section symmetries of the problem are described. These symmetries help us to find
some Maxwell points.

{\Def
\label{def-symmetry}
\emph{A symmetry} of the exponential map is a pair of diffeomorphisms
$$
s: C \times \R_+ \rightarrow C \times \R_+ \qquad
\text{and} \qquad
\widehat{s} : G \rightarrow G \qquad
\text{such that} \qquad
\Exp \circ s = \widehat{s} \circ \Exp.
$$
}

Next we consider only symmetries that correspond to isometries of $C$ (in the sense of the Killing form) that conserve or invert the vertical part of the Hamiltonian vector field
$$
\vec{H}_{vert} (p) = - \frac{\eta p_3}{I_1} (-p_2 \varepsilon_1 + p_1 \varepsilon_2) \in \VectField(C).
$$
It is clear that the group of such isometries is generated by rotations around the axis $\sspan\{e_3\}$ and the reflections $\sigma_1$ and $\sigma_2$ in the planes
$\sspan\{e_1, e_3\}$ and $\sspan\{e_1, e_2\}$ respectively. Denote this group by $S$. It is isomorphic to $\mathrm{O}_2 \times \mathbb{Z}_2$.

Let us introduce the following notation:
\begin{equation*}
\begin{array}{rcl}
Q^e_{t, p} & = & \cos{t} + \sin{t} (p_1 i + p_2 j + p_3 k), \\
Q^p_{t, p} & = & 1 + \frac{t}{2} (p_1 i + p_2 j + p_3 k), \\
Q^h_{t, p} & = & \cosh{t} + \sinh{t} (p_1 i + p_2 j + p_3 k). \\
\end{array}
\end{equation*}
Note that any element of $\SU_{1, 1}$ has one of the following forms:
$Q^e_{t, p}, \pm Q^p_{t, p}, \pm Q^h_{t, p}$, where $t \geqslant 0$.
The parameter $t$ is unique up to addition of $2\pi$ in the case of $Q_{t, p}^e$.

{\Proposition
\label{prop-symmeties}
The group $S$ is embedded into the group of symmetries of the exponential map. To any element $\sigma \in S$ assign the pair of diffeomorphisms
$$
s_{\sigma} : C \times \mathbb{R}_{+} \rightarrow C \times \mathbb{R}_{+} \qquad \text{and} \qquad
\widehat{s}_{\sigma} : G \rightarrow G
$$
given by
\begin{equation*}
s_{\sigma}(p, t) =
\left\{
\begin{array}{lll}
(\sigma(p), t), &  \text{for} & d\sigma(\vec{H}_{vert}) = \vec{H}_{vert}, \\
(\sigma e^{t \vec{H}_{vert}} (p), t), & \text{for} & d\sigma(\vec{H}_{vert}) = -\vec{H}_{vert},
\end{array}
\right.
\end{equation*}
$$
\widehat{s}_{\sigma}(\Pi Q^e_{t, p}) = \Pi Q^e_{t, \sigma (p)}, \qquad
\widehat{s}_{\sigma}(\Pi Q^p_{t, p}) = \Pi Q^p_{t, \sigma (p)}, \qquad
\widehat{s}_{\sigma}(\Pi Q^h_{t, p}) = \Pi Q^h_{t, \sigma (p)}.
$$
}
\medskip

\emph{Proof.}
Note that the action of $\widehat{s}_{\sigma}$ does not depend on the choice of pre-image of $s \in \PSL_2(\R)$ under the covering $\Pi$.

It is enough to check that for generators $\sigma \in S$ the pair of diffeomorphisms $(s_{\sigma}, \widehat{s}_{\sigma})$ is a symmetry of the exponential map.
Generators of $S$ are rotations around the line $\sspan\{e_3\}$ and reflections in the planes
$\sspan\{e_1, e_3\}$ and $\sspan\{e_1, e_2\}$.

Such rotations and the first reflection do not change $p_3$,
therefore they do not change the components $q_0$ and $q_3$ of a corresponding split-quaternion. The second reflection changes the sign of $p_3$,
so the component $q_0$ does not change, but the component $q_3$ changes the sign.

Thus, we need to know how the components $q_1$ and $q_2$ of the endpoint of the geodesic change
when the generators of $S$ act on the initial momentum of the geodesic. It is enough to show that
$$
\left(
\begin{array}{c}
q_1 (t, \sigma (p) ) \\
q_2 (t, \sigma (p) ) \\
\end{array}
\right)
=
\sigma
\left(
\begin{array}{c}
q_1 (t, p) \\
q_2 (t, p) \\
\end{array}
\right).
$$
From~(\ref{eq-time-geodesic}, \ref{eq-space-geodesic}, \ref{eq-light-geodesic}) one can see that this is true for rotations around the axis $\sspan\{e_3\}$, since the transformation $R_{e_3, \frac{t \eta p_3}{2 I_1}}$ is such a rotation and it commutes with $\sigma$.

If $\sigma$ is one of reflections then it reverses the vertical part of the Hamiltonian vector field. Hence
$$
\Exp \circ s_{\sigma}(p, t) = \Exp (\sigma e^{t \vec{H}_{vert}}(p), t) =
\Exp (\sigma R_{e_3, -\frac{t\eta p_3}{I_1}}, t).
$$

If $\sigma$ is reflection in the plane $\sspan\{e_1, e_3\}$ then
$R_{e_3, \varphi} \sigma = \sigma R_{e_3, -\varphi}$ and $\sigma(p)_3 = p_3$, thus
$$
\left(
\begin{array}{c}
q_1 (t, \sigma (p) ) \\
q_2 (t, \sigma (p) ) \\
\end{array}
\right)
=
R_{e_3, -\frac{t \eta p_3}{2 I_1}}
\sigma R_{e_3, -\frac{t\eta p_3}{I_1}}
\left(
\begin{array}{c}
p_1 \\
p_2 \\
\end{array}
\right)
=
\sigma R_{e_3, -\frac{t \eta p_3}{2 I_1}}
\left(
\begin{array}{c}
p_1 \\
p_2 \\
\end{array}
\right)
=
\sigma
\left(
\begin{array}{c}
q_1 (t, p) \\
q_2 (t, p) \\
\end{array}
\right).
$$

If $\sigma$ is reflection in the plane $\sspan\{e_1, e_2\}$ then
$\sigma$ commutes with rotations around the axis $\sspan\{e_3\}$, but $\sigma(p)_3 = -p_3$, whence
$$
\left(
\begin{array}{c}
q_1 (t, \sigma (p) ) \\
q_2 (t, \sigma (p) ) \\
\end{array}
\right)
=
R_{e_3, -\frac{t \eta (-p_3)}{2 I_1}}
\sigma
R_{e_3, - \frac{t\eta p_3}{I_1}}
\left(
\begin{array}{c}
p_1 \\
p_2 \\
\end{array}
\right)
=
\sigma R_{e_3, -\frac{t \eta p_3}{2 I_1}}
\left(
\begin{array}{c}
p_1 \\
p_2 \\
\end{array}
\right)
=
\sigma
\left(
\begin{array}{c}
q_1 (t, p) \\
q_2 (t, p) \\
\end{array}
\right).
$$

Hereby we have shown that $\Exp \circ s_{\sigma}(p, t) = \widehat{s}_{\sigma} \circ \Exp (p, t)$ for generators $\sigma \in S$.
\qquad$\Box$

\section{\label{section-maxwell-strata}Maxwell strata}

{\Def
\label{def-maxvell-point-and-time}
\emph{A Maxwell point} is a point $Q \in G$ such that there are two distinct geodesics with arclength parametrization $Q_1, Q_2: [0, T] \rightarrow \SO_3$, coming to the point $Q$ at the same time $Q = Q_1(t_{max}) = Q_2(t_{max})$. This time is called \emph{a Maxwell time}.}

It is known (see for example~\cite{sachkov-didona1}) that after a Maxwell point an extremal trajectory can not be optimal.

{\Def
\label{def-first-maxwell-point-and-time}
\emph{The first Maxwell set in the pre-image of the exponential map} is the set
\begin{multline*}
\M = \{(p, t_{max}) \in C \times \R_+ \ | \ \exists p' \in C \setminus \{p\} :
\Exp(p, t_{max}) = \Exp(p', t_{max}), \\
\text{but} \ \forall t \in (0, t_{max}) \ \forall p_1 \in C \setminus \{p\} \
\Exp(p, t) \neq \Exp(p_1, t)\}.
\end{multline*}
The time $t_{max}$ is called \emph{the first Maxwell time for} $p \in C$.}

Obviously, $\Exp{(\M)}$ consists of Maxwell points.

{\Def
\label{def-first-maxwell-point-and-time-symmetries}
Suppose $A$ is a subset of the group $S$. \emph{The first Maxwell set that corresponds to $A$ in the pre-image of the exponential map} is the set
\begin{multline*}
\M(A) = \{(p, t_{max}) \in C \times \R_+ \ | \ \exists \sigma \in A \ : \
\Exp(p, t_{max}) = \Exp \circ s_{\sigma} (p, t_{max}), \\
\text{but} \ \forall t \in (0, t_{max}) \
\forall \sigma \in A\setminus\{\id\} \
\Exp(p, t) \neq \Exp \circ s_{\sigma} (p, t) \}.
\end{multline*}
The time $t_{max}$ is called \emph{the first Maxwell time corresponding to $A$ for} $p \in C$.}

This time is not less than the first Maxwell time.

The aim of this section is description of the first Maxwell strata in the image and pre-image of the exponential map. First we describe the sets $\M(\sigma)$ for each $\sigma \in S$.
Second we explore the relative location of the sets $\M(\sigma)$ and then find
$$
\M(S) \subset \bigcup_{\sigma \in S} \M(\sigma).
$$

Next we will show that the exponential map is a diffeomorphism from the domain of $C \times \R_+$
bounded by $\overline{\M(S)}$ to $G \setminus (\Exp{\overline{\M(S)}} \cup \{\id\})$. This will imply that $\overline{\M(S)}$ and $\Exp{\overline{\M(S)}}$ are the cut loci in the pre-image and image of the exponential map respectively. This means that
$\overline{\M(S)} = \{(p, t_{cut}(p)) \in C \times \R_+ \}$,  where $t_{cut}(p)$ is a time such that the geodesic $\{\Exp(p, t) \ | \ t \in \R_+\}$ is a shortest arc for $t \in [0, t_{cut}(p)]$ but it is not a shortest arc for $t > t_{cut}(p)$.

\subsection{\label{section-maxwell-strata-corresponding-to-symmetries}Maxwell strata corresponding to symmetries}

{\Def
\label{def-first-positive-root}
Denote by $C^e, C^p, C^h$ the subsets of $C$ consisting of time-, light- or space-like covectors respectively.
Introduce the following notation:
$$
\begin{array}{lclclcl}
\tau^e_0(\p) & = & \min\{\tau \in \R_+ \ | \ q_0^e(\tau, \p) = 0 \}, & \qquad &
\tau^e_3(\p) & = & \min\{\tau \in \R_+ \ | \ q_3^e(\tau, \p) = 0 \},\\
t^p_0(p) & = & \min\{t \in \R_+ \ | \ q_0^p(t, p) = 0 \}, & \qquad &
t^p_3(p) & = & \min\{t \in \R_+ \ | \ q_3^p(t, p) = 0 \},\\
\tau^h_0(\p) & = & \min\{\tau \in \R_+ \ | \ q_0^h(\tau, \p) = 0 \}, & \qquad &
\tau^h_3(\p) & = & \min\{\tau \in \R_+ \ | \ q_3^h(\tau, \p) = 0 \}.
\end{array}
$$
}

We consider the values of $\tau^e_0, \tau^e_3, \tau^h_0, \tau^h_3$ as functions of the variable $\p$ and the parameter $\eta$. For time- and space-like covectors $p$ the functions $q_0(\tau)$ and $q_3(\tau)$ are even and odd respectively. Thus, all of the functions $\tau^e_0, \tau^e_3, \tau^h_0, \tau^h_3$ are even. For $\p = 0$ the equations $q^h_0(\tau, \p) = 0$ and $q^h_3(\tau, \p) = 0$
read as $\cosh \tau = 0$ and identity respectively. That is why the values $\tau^h_0(0)$ and $\tau^h_3(0)$ are undefined. Therefore, we can consider the following domains of
the functions:
\begin{equation*}
\begin{array}{rcl}
\tau^e_0: [1, +\infty) \rightarrow \R_+, & \qquad & \tau^h_0: (0, +\infty) \rightarrow \R_+, \\
\tau^e_3: [1, +\infty) \rightarrow \R_+, & \qquad & \tau^h_3: (0, +\infty) \rightarrow \R_+. \\
\end{array}
\end{equation*}

{\Proposition
\label{prop-symmetries-maxwell-strata}
The set
$$
\bigcup_{\sigma \in S} \M(\sigma) = \M_0 \cup \M_{12} \cup \M_3
$$
is the union of the three strata
$$
\begin{array}{lcl}
\M_0 & = & \M^e_0 \cup \M^p_0 \cup \M^h_0, \\
\M_{12} & = & \{ (p, t) \in C^e \times \R_+ \ | \ \p \neq \pm 1, \ t = \frac{2 \pi I_1}{|p|} \}, \\
\M_3 & = & \M^e_3 \cup \M^p_3 \cup \M^h_3,
\end{array}
$$
where
$$
\begin{array}{ll}
\M^e_0  =  \{ (p, t) \in C^e \times \R_+ \ | \ t = \frac{2 \tau^e_0(|\p|) I_1}{|p|} \}, &
\M^e_3  =  \{ (p, t) \in C^e \times \R_+ \ | \ t = \frac{2 \tau^e_3(|\p|) I_1}{|p|} \}, \\
\M^p_0  =  \{ (p, t) \in C^p \times \R_+ \ | \ t = t^p_0(p) \}, &
\M^p_3  =  \{ (p, t) \in C^p \times \R_+ \ | \ t = t^p_3(p) \}, \\
\M^h_0  =  \{ (p, t) \in C^h \times \R_+ \ | \ \p \neq 0, &
\M^h_3  =  \{ (p, t) \in C^h \times \R_+ \ | \ \p \neq 0, \\
\qquad \qquad \qquad \qquad \qquad \qquad t = \frac{2 \tau^h_0(|\p|) I_1}{|p|} \}, &
\qquad \qquad \qquad \qquad \qquad \qquad t = \frac{2 \tau^h_3(|\p|) I_1}{|p|} \}. \\
\end{array}
$$
}
\medskip

\emph{Proof.}
It is clear that $\M(\sigma) \subset G^{\sigma} = \{ g \in G \ | \ \widehat{s}_{\sigma} (g) = g \}$. For any $\sigma \in S$ consider the set of its fixed points $G^{\sigma}$. For which of them are there  two symmetric geodesics coming there at the same time?

Evidently the sets of fixed points in $G$ for different elements of $S$ lie in the union of the sets
$$
\Pi \{q \in \SU_{1, 1} \ | \ q_0 = 0 \}, \qquad
\Pi \{q \in \SU_{1, 1} \ | \ q_1 = q_2 = 0 \}, \qquad
\Pi \{q \in \SU_{1, 1} \ | \ q_3 = 0 \}.
$$
For any covector in $C$ there is a symmetric (with respect to some element of $S$) covector such that the two geodesics with these initial momenta come to one of these sets at the same time. This time is equal to the first positive root of the corresponding equation. The covectors with $\p = 0$ are exceptions: a geodesic with such initial momentum always lies in
$\Pi \{q \in \SU_{1, 1} \ | \ q_3 = 0 \}$ and never reaches $\Pi \{q \in \SU_{1, 1} \ | \ q_0 = 0 \}$. Note that only geodesics
with initial momenta from $C^e$ reach the set
$\Pi \{q \in \SU_{1, 1} \ | \ q_1 = q_2 = 0 \}$, and the geodesics with $\p = \pm 1$ always lie in this set. For details see the similar Proposition~2 in paper~\cite{podobryaev-sachkov} about the Riemannian problem on $\SO_3$.
\qquad$\Box$
\medskip

\subsection{\label{section-some-propertis-of-maxwell-time}The functions
$\tau^e_0, \tau^e_3, \tau^h_0, \tau^h_3$ are continuous}

To investigate the relative location of the Maxwell strata $\M_0$, $\M_{12}$ and $\M_3$ we need to compare the corresponding Maxwell times: for time-like initial momenta $\frac{2 \tau^e_0(\p) I_1}{|p|}$, $\frac{2 \pi I_1}{|p|}$ and $\frac{2 \tau^e_3(\p) I_1}{|p|}$; for light-like initial momenta $t^p_0(p), t^p_3(p)$; and for space-like ones $\frac{2 \tau^h_0(\p) I_1}{|p|}, \frac{2 \tau^h_3(\p) I_1}{|p|}$. Since $|p|$ depends only on $\p$, it is enough to compare the functions
$\tau^e_0(\p), \tau^e_3(\p)$ and the number $\pi$ for different values of $\p \in [1, +\infty)$, and compare $\tau^h_0(\p)$ and $\tau^h_3(\p)$ for $\p \in (0, +\infty)$. For this purpose let us examine some properties of these functions.

{\Proposition
\label{prop-continious}
The functions $\tau^e_0, \tau^e_3, \tau^h_0, \tau^h_3$ are continuous on their domains.
}
\medskip

\emph{Proof.}
The implicit function theorem implies that
it is enough to verify that the functions $q_0$ and $q_3$ have no multiple roots. (We consider $q_0$ and $q_3$ as functions of the variable $\tau$ and the parameter $\p$.) Let us check this for time- and space-like parameters $p$ together. Introduce some notation to make computations more easy:
\begin{equation}
\label{def-trig-hyp}
c(\tau, p) =
\begin{aligned}
\left\{
\begin{array}{lll}
\cos{\tau}, & \text{for} & \type{(p)} = 1, \\
\cosh{\tau}, & \text{for} & \type{(p)} = -1, \\
\end{array}
\right.
\end{aligned}
\qquad
s(\tau, p) =
\begin{aligned}
\left\{
\begin{array}{lll}
\sin{\tau}, & \text{for} & \type{(p)} = 1, \\
\sinh{\tau}, & \text{for} & \type{(p)} = -1. \\
\end{array}
\right.
\end{aligned}
\end{equation}
Then the following equations hold:
$$
c^2(\tau, p) + \type(p) s^2(\tau, p) = 1,
$$
$$
\frac{\partial c}{\partial \tau}(\tau, p) = -\type(p) s(\tau, p), \qquad
\frac{\partial s}{\partial \tau}(\tau, p) = c(\tau, p),
$$
where $\type(p) = \sgn(-\Kil(p))$.

Calculate derivatives of the functions $q_0$ and $q_3$ of the variable $\tau$:
$$
\frac{\partial q_0}{\partial \tau} =
-\type(p) (1 + \type(p) \eta \p^2) s(\tau, p) \ce - \p (1 + \eta) c(\tau, p) \se,
$$
$$
\frac{\partial q_3}{\partial \tau} =
-\type(p)(1 + \type
(p) \eta \p^2) s(\tau, p) \se + \p (1 + \eta) c(\tau, p) \ce.
$$

1. Assume that $q_0(\tau)$ has a multiple root. This means that there is $\p \in [0, +\infty)$ such that
\begin{equation}
\label{eq-multiroot-q0}
\left\{
\begin{aligned}
q_0(\tau) = 0, \\
\frac{\partial q_0}{\partial \tau}(\tau) = 0.
\end{aligned}
\right.
\end{equation}
Let us divide the both equations by $c(\tau, p) \cos (\tau \eta \p)$ and denote $t(\tau, p) = \frac{s(\tau, p)}{c(\tau, p)}$. The case when the denominator equals zero will be considered below. We have
\begin{equation*}
\left\{
\begin{aligned}
1 - \p t(\tau, p) \tan (\tau \eta \p) = 0,\\
-\type(p)(1 + \type(p) \eta \p^2) t(\tau, p) - \p (1 + \eta) \tan (\tau \eta \p) = 0.
\end{aligned}
\right.
\end{equation*}
Note that $t(\tau, p) \neq 0$ and $\p \neq 0$. Expressing $\tan (\tau \eta \p)$ from the first equation and substituting it to the second one, we obtain
$$
t^2 (\tau, p) = - \frac{1 + \eta}{\type(p)(1 + \type(p) \eta \p^2)}.
$$
Let us show that $\type(p)(1 + \type(p) \eta \p^2) < 0$. Indeed, if $\type(p) > 0$ then $\p \geqslant 1$, and  $\eta < -1$ implies $\eta \p^2 < -1$. When $\type(p) < 0$ we have $\p > 0$ and $-\eta \p^2 > 0$, thus $(1 + \type(p) \eta \p^2) > 0$. Besides  $1 + \eta < 0$. Hence $t^2 (\tau, p) < 0$. We get a contradiction.

Now consider the case when the denominator $c(\tau, p) \cos (\tau \eta \p)$ equals zero.
If $c(\tau, p) = 0$, then from system~(\ref{eq-multiroot-q0}) we have
\begin{equation*}
\left\{
\begin{aligned}
-\p \se = 0, \\
-\type(p) (1 + \type(p) \eta \p^2) \ce = 0.
\end{aligned}
\right.
\end{equation*}
Since cosine and sine can not be zero simultaneously and $\type(p) (1 + \type(p) \eta \p^2) < 0$, we obtain $\p = 0$, thus $\ce = 1$, and we get a contradiction.

If $\ce = 0$, then from system~(\ref{eq-multiroot-q0}) we get
\begin{equation*}
\left\{
\begin{aligned}
-\p s (\tau, p) = 0, \\
-\p (1 + \eta) c (\tau, p) = 0,
\end{aligned}
\right.
\end{equation*}
thus $\p = 0$. Then $\ce
 = 1$, we get a contradiction.

2. Assume that $q_3$ has a multiple root. Thus, for some $\p \in (0, +\infty)$ we have
\begin{equation}
\label{eq-multiroot-q3}
\left\{
\begin{aligned}
q_3(\tau) = 0, \\
\frac{\partial q_3}{\partial \tau}(\tau) = 0.
\end{aligned}
\right.
\end{equation}
If $c(\tau, p) \cos (\tau \eta \p)$ is non zero, then divide both equations by this expression. We get
\begin{equation*}
\left\{
\begin{aligned}
\tan (\tau \eta \p) + \p t(\tau, p) = 0,\\
-\type(p)(1 + \type(p) \eta \p^2) t(\tau, p) \tan (\tau \eta \p) + \p (1 + \eta) = 0.
\end{aligned}
\right.
\end{equation*}
Since $\p \neq 0$ we have
$$
t^2 (\tau, p) = - \frac{1 + \eta}{\type(p)(1 + \type(p) \eta \p^2)},
$$
this fraction is less than zero (see item~1), we get a contradiction.

Consider now the case when the denominator $c(\tau, p) \cos (\tau \eta \p)$ is equal to zero.

If $c(\tau, p) = 0$, then from system~(\ref{eq-multiroot-q3}) we get
\begin{equation*}
\left\{
\begin{aligned}
-\p \ce = 0, \\
-\type(p) (1 + \type(p) \eta \p^2) \se = 0,
\end{aligned}
\right.
\end{equation*}
hence (since cosine and sine can not be equal to zero simultaneously and $\type(p) (1 + \type(p) \eta \p^2) < 0$) we have $\p = 0$, a contradiction.

If $\ce = 0$, then from system~(\ref{eq-multiroot-q3}) we have
\begin{equation*}
\left\{
\begin{aligned}
c (\tau, p) = 0, \\
-\type(p) (1 + \type(p) \eta \p^2) s (\tau, p) = 0,
\end{aligned}
\right.
\end{equation*}
we get a contradiction.
\qquad$\Box$
\medskip

\subsection{\label{section-maxwell-strata-location}Relative location of Maxwell strata}

Now we compare $\tau^e_0(\p)$, $\pi$ and $\tau^e_3(\p)$ for different values of
$\p \in [1, +\infty)$ and compare $\tau^h_0(\p)$ and $\tau^h_3(\p)$ for
$\p \in (0, +\infty)$. Thereby we will explore the relative location of the Maxwell strata.

{\Proposition
\label{prop-comparsion-e-tau0-tau3}
For all $\p \in [1, +\infty)$ the inequality $\tau^e_0(\p) < \tau^e_3(\p)$ is satisfied.
}
\medskip

\emph{Proof.}
Notice that for $\p = 1$ the statement of the proposition is true. Indeed,
$$
q_0^e(\tau) = \cos(\tau (1 + \eta)), \qquad
q_3^e(\tau) = \sin(\tau (1 + \eta)).
$$
Then $\tau^e_0(1) = -\frac{\pi}{2 (1 + \eta)} < -\frac{\pi}{(1 + \eta)} = \tau^e_3(1)$.

Assume (by contradiction) that for some $\p$ there holds the inequality
$\tau^e_0(\p) \geqslant \tau^e_3(\p)$. Because of continuity of the functions $\tau^e_0$ and $\tau^e_3$ (Proposition~\ref{prop-continious}) there is $\widehat{p}_3$ such that $\tau^e_0(\widehat{p}_3) = \tau^e_3(\widehat{p}_3)$. This means that for some $\widehat{p}_3$ and $\tau$ we have
$q_0(\tau, \widehat{p}) = q_3(\tau, \widehat{p}) = 0$ in contradiction with $q_0^2 - q_1^2 - q_2^2 + q_3^2 = 1$.
\qquad$\Box$
\medskip

The above proposition shows that geodesics with time-like initial momenta reach the stratum $\Exp{(\M^e_0)}$ earlier than the stratum $\Exp{(\M^e_3)}$.

Consider now the strata $\M^e_0$ and $\M^e_{12}$.

\begin{figure}[h]
\caption{The function $\tau^e_0(\p)$ and $\pi$.}
\label{pic-maxwell-time}
\medskip
     \begin{minipage}[h]{0.30\linewidth}
        \center{\includegraphics[width=1\linewidth]{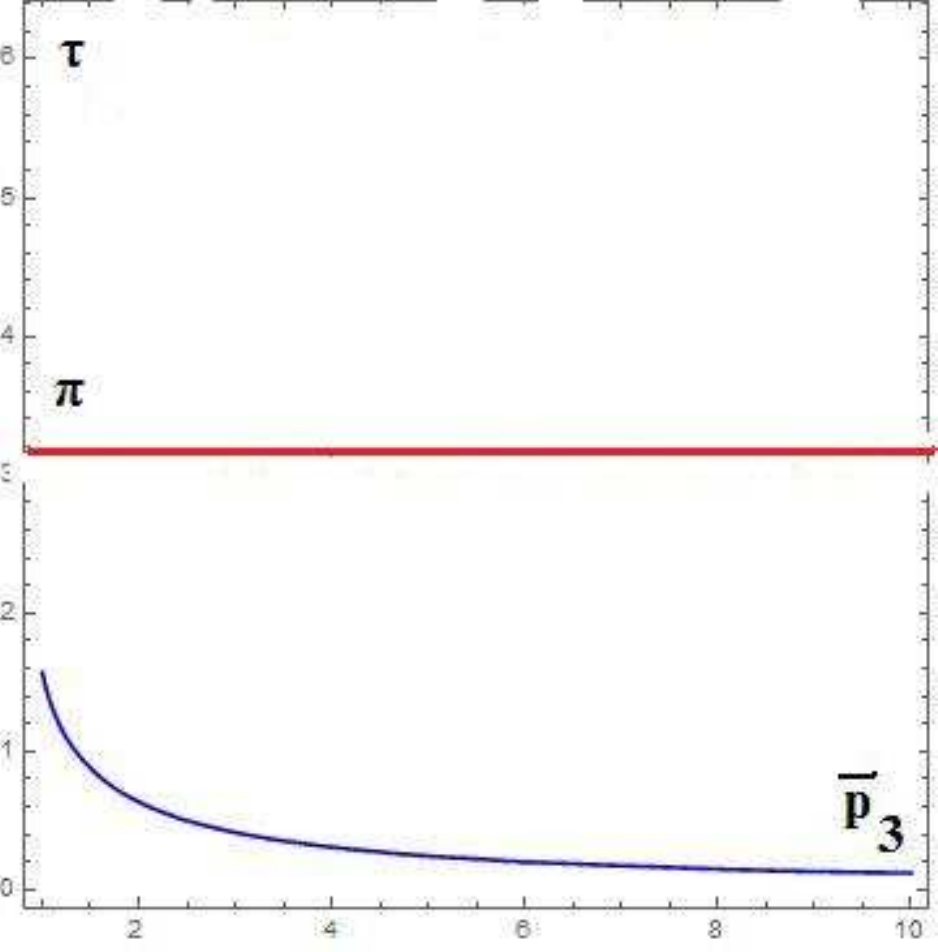} \\ $\eta < -\frac{3}{2}$}
     \end{minipage}
     \hfill
     \begin{minipage}[h]{0.30\linewidth}
        \center{\includegraphics[width=1\linewidth]{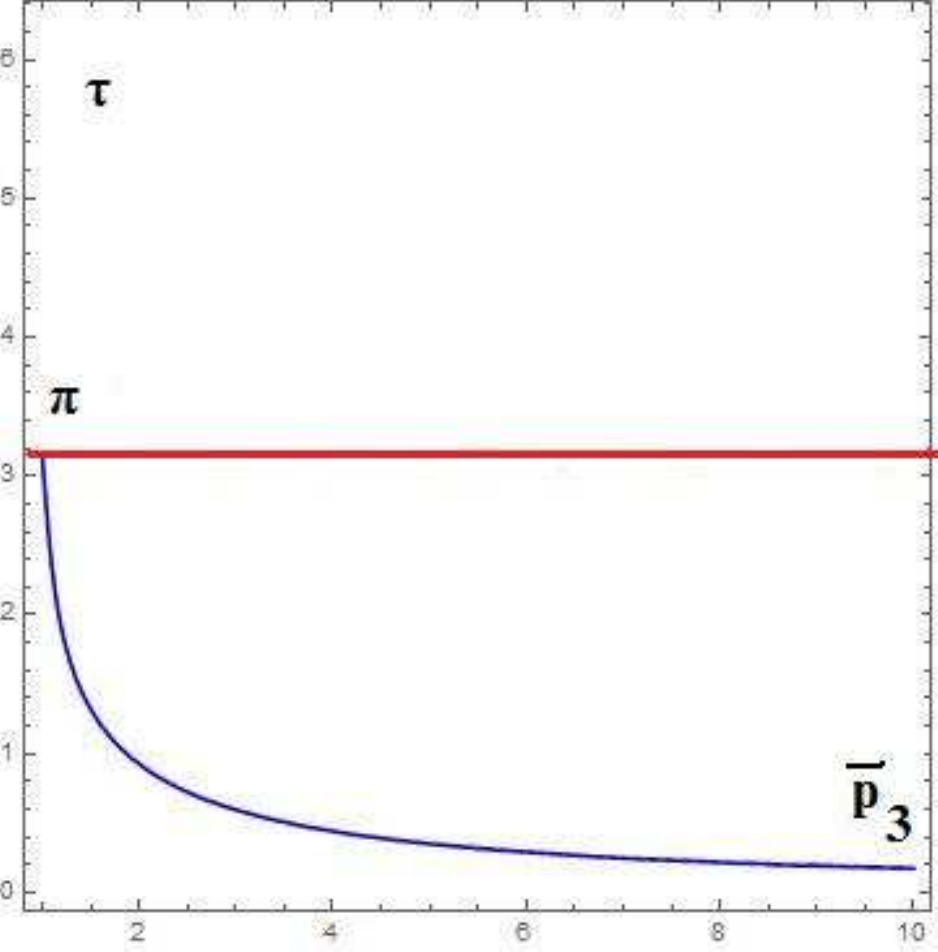} \\ $\eta = -\frac{3}{2}$}
     \end{minipage}
     \hfill
     \begin{minipage}[h]{0.30\linewidth}
        \center{\includegraphics[width=1\linewidth]{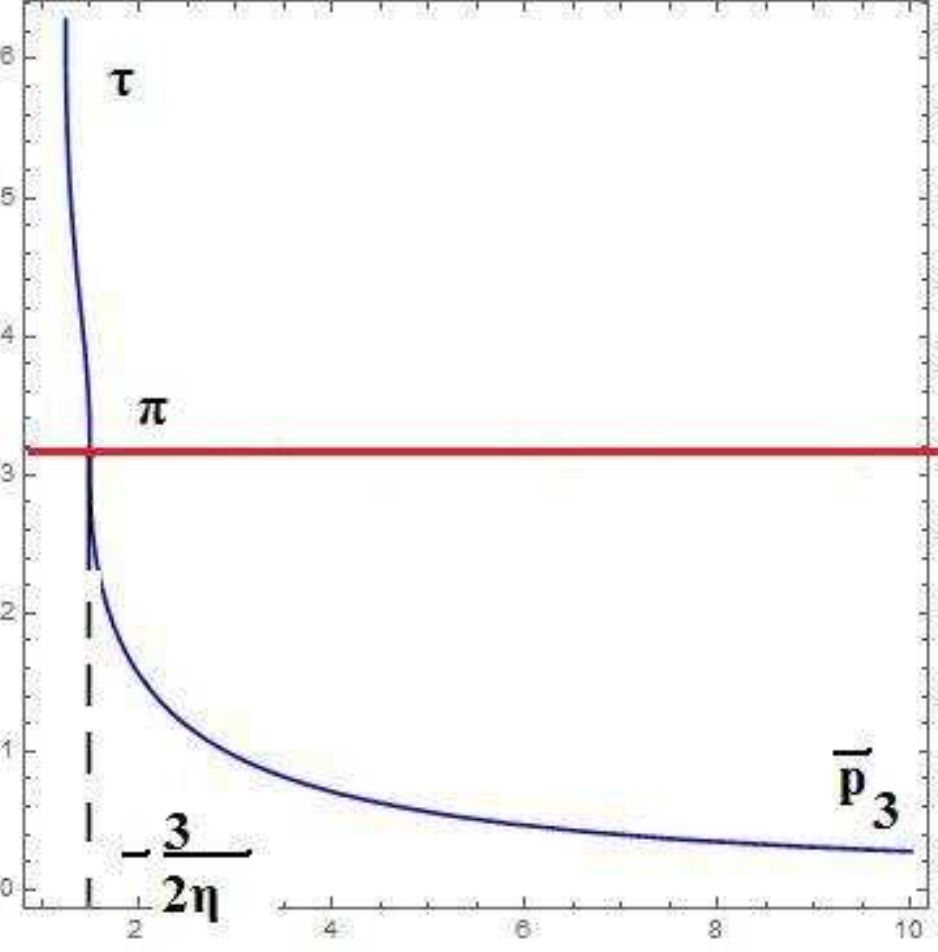} \\ $-\frac{3}{2} < \eta $}
     \end{minipage}
\end{figure}

{\Proposition
\label{prop-comparsion-e-tau0-pi}
$(1)$ If $\eta \leqslant -\frac{3}{2}$, then for all $\p \in [1, +\infty)$ the inequality $\tau^e_0(\p) \leqslant \pi$ is satisfied.\\
$(2)$ If $\eta > -\frac{3}{2}$, then
$\tau^e_0(\p) \geqslant \pi$ for $\p \in [0, -\frac{3}{2 \eta}]$ and
$\tau^e_0(\p) < \pi$ for $\p \in (-\frac{3}{2 \eta}, +\infty)$.\\
See Figure~$\ref{pic-maxwell-time}$.
}
\medskip

\emph{Proof.}
(1) Note that $q_0^e(0) = 1$. Hence, it is enough to find $\theta \in (0, \pi]$ such that $q_0^e(\theta) \leqslant 0$, since in this case the continuous function $q_0^e$ of the variable $\tau$ has zero at the interval $(0, \theta]$, i.e., a zero that is less than or equal to $\pi$. Take
\begin{equation*}
\theta =
\left\{
\begin{array}{rl}
\pi, & \text{for} -\eta \p < 2, \\
-\frac{\pi}{\eta \p}, & \text{for} -\eta \p \geqslant 2.
\end{array}
\right.
\end{equation*}
When $\eta \leqslant -\frac{3}{2}$, for $\p \geqslant 1$ we have $-\eta \p \geqslant \frac{3}{2}$. Then $-\frac{\pi}{\eta \p} \leqslant \frac{2 \pi}{3} < \pi$. Thus
\begin{equation*}
q_0^e(\theta) =
\left\{
\begin{array}{rl}
-\cos(\pi \eta \p), & \text{for} -\eta \p < 2, \\
-\cos(-\frac{\pi}{\eta \p}), & \text{for} -\eta \p \geqslant 2.
\end{array}
\right.
\end{equation*}
In the first case $-2 \pi < \pi \eta \p \leqslant -\frac{3 \pi}{2}$, then $\cos(\pi \eta \p) \geqslant 0$. In the second case $0 < -\frac{\pi}{\eta \p} \leqslant \frac{\pi}{2}$, then $\cos(-\frac{\pi}{\eta\p}) \geqslant 0$. So, we get $q_0^e(\theta) \leqslant 0$.

(2) Firstly obtain the second part of the statement. For $\p \geqslant -\frac{3}{2 \eta}$ we have $-\eta \p \geqslant \frac{3}{2}$. Choose now the same $\theta$ as in the proof of statement~(1).

Obtain now the first part of statement~(2). When $\p = 1$ we get $q_0^e(\tau) = \cos(\tau(1+\eta))$. The first positive root of this function is $-\frac{\pi}{2 (1+\eta)} > \pi$. Assume (by contradiction) that there exists $\p' \in (1, -\frac{3}{2 \eta})$ such that $\tau_0^e(\p') < \pi$. Because of continuity of the function $\tau_0^e$ there is $\widehat{p}_3 \in (1, \p')$ such that $q_0^e(\pi, \widehat{p}_3) = -\cos(\pi \eta \widehat{p}_3) = 0$. Whence $\widehat{p}_3 = \frac{2k+1}{2\eta}, \ k \in \Z$. It is easy to see that for all $k \in \Z$ the point $\widehat{p}_3$ is outside of the interval $(1, -\frac{3}{2 \eta})$. We get a contradiction.
\qquad$\Box$
\medskip

{\Proposition
\label{prop-comparsion-h-tau0-tau3}
For all $\p \in (0, +\infty)$ the inequality $\tau_0^h(\p) < \tau_3^h(\p)$ is satisfied.
}
\medskip

\emph{Proof.}
Actually we need to find at least one $\p$ such that $\tau_0^h(\p) < \tau_3^h(\p)$.
Indeed, if there is a point  such that this inequality is violated, then there exists a point at which the continuous functions $\tau_0^h$ and $\tau_3^h$ have the same value. This means that $q_0$ and $q_3$ vanish simultaneously in contradiction with $q_0^2 - q_1^2 - q_2^2 + q_3^2 = 1$.

Let us verify that the required inequality holds for $\p = 1$. We have
$$
q_0^h(\tau_0^h(1), 1) = \cosh \tau_0^h(1) \cos(\tau_0^h(1)\eta) - \sinh \tau_0^h(1) \sin(\tau_0^h(1)\eta) = 0,
$$
$$
q_3^h(\tau_0^h(1), 1) = \cosh \tau_3^h(1) \sin(\tau_3^h(1)\eta) + \sinh \tau_3^h(1) \cos(\tau_3^h(1)\eta) = 0.
$$
Note that $\sin(\tau_0^h(1)\eta) \neq 0$, since otherwise from the first equation we have $\cos(\tau_0^h(1)\eta) = 0$, a contradiction. Furthermore $\cos(\tau_3^h(1)\eta) \neq 0$, since otherwise from the second equation we have $\sin(\tau_3^h(1)\eta) = 0$, a contradiction. The hyperbolic cosine never vanishes, so the expressions $\cosh \tau_0^h(1) \sin (\tau_0^h(1)\eta)$ and $\cosh \tau_3^h(1) \cos (\tau_3^h(1)\eta)$ never vanish as well. Divide the first and the second equations by these expressions respectively. Now we need to compare the first positive roots of the equations
$$
\cot (\tau\eta) = \tanh \tau,
$$
$$
-\tan (\tau \eta) = \tanh \tau.
$$
The first positive root of the first equation lies inside of the interval $(-\frac{\pi}{2\eta}, -\frac{\pi}{\eta})$. The first positive root of the second one lies in the interval $(-\frac{\pi}{\eta}, -\frac{3\pi}{2\eta})$. (We use the fact that the derivative of the function $-\tan (\tau \eta)$ at zero is equal to $-\eta > 1$, i.e., it is greater than the derivative of the function $\tanh \tau$ at zero. Therefore, the graph of the function $\tanh \tau$ intersects the first branch of the graph of the function $-\tan (\tau \eta)$ only at zero.) So, we have $\tau_0^h(1) < \tau_3^h(1)$.
\qquad$\Box$
\medskip

The above proposition implies that geodesics with space-like initial momenta reach the stratum $\Exp{(\M^h_0)}$ earlier than the stratum $\Exp{(\M^h_3)}$. The proposition below states that the same is true for light-like initial momenta.

{\Proposition
\label{prop-comparsion-p-t0-t3}
For $p \in C^p$ the inequality $t^p_0(p) < t^p_3(p)$ holds.
}
\medskip

\emph{Proof.}
The equation $q_0^p(t) = 0$ reads as
$$
\ctl - \frac{t p_3}{2 I_1} \stl = 0.
$$
It is easy to see that $\sin (\frac{t \eta p_3}{2 I_1}) \neq 0$. Divide the equation by this expression. Denote $\tau_p = \frac{t p_3}{2 I_1}$. We get an equivalent equation
$$
\cot(\tau_p \eta) = \tau_p.
$$
Its first positive root lies in the interval $(-\frac{\pi}{2\eta}, -\frac{\pi}{\eta})$.

Similarly, the equation $q_3^p = 0$ is equivalent to the equation
$$
\tan(\tau_p \eta) = -\tau_p.
$$
Its first positive root lies in $(-\frac{2\pi}{\eta}, -\frac{5\pi}{2\eta})$.
(Since derivative of the function $\tan (\tau_p \eta)$ at zero is equal to $-\eta > 1$, hence
the graph of the function $-\tau_p$ intersects the first branch of the graph of the function $\tan (\tau_p \eta)$ at zero only.)

The last interval is located to the right of the first one, so we get the statement of the proposition.
\qquad$\Box$
\medskip

Denote by $t_{max}(p), \ p \in C,$ the first Maxwell time corresponding to the symmetry group $S$ of the exponential map. Propositions~\ref{prop-comparsion-e-tau0-tau3}, \ref{prop-comparsion-e-tau0-pi}, \ref{prop-comparsion-h-tau0-tau3}, \ref{prop-comparsion-p-t0-t3} imply that
\begin{equation}
\label{eq-maxwell-time}
t_{max} (p) = \left \{
\begin{array}{rcll}
\frac{2 I_1}{|p|}\min{(\pi, \tau^e_0(|\p|))}, & \text{for} & p \in C^e, &  \\
t^p_0(\p), & \text{for} & p \in C^p, & \\
\frac{2 I_1}{|p|} \tau^h_0(|\p|), & \text{for} & p \in C^h, & \p \neq 0,\\
+\infty, & \text{for} & p \in C^h, & \p = 0.\\
\end{array}
\right.
\end{equation}

{\Lemma
\label{lemma-firt-maxwell-time-continious}
The function $t_{max} : C \rightarrow (0, +\infty]$ is continuous.
}
\medskip

\emph{Proof.}
Proposition~\ref{prop-continious} implies that it is enough to proof that:
\begin{enumerate}
\item the first Maxwell time is continuous at points $p \in C$ with $|p| = 0$;
\item there holds the equality: $\lim_{\p \rightarrow 0}{\frac{2 \tau^h_0(|\p|) I_1}{|p|}} = +\infty$.
\end{enumerate}

1. The map $\Exp : C \times \R_+ \rightarrow \SU_{1, 1}$ is smooth (Remark~\ref{rem-exponential-map-is-real-analitic}), hence its component $q_0$ (in the coordinates $q_0, q_1, q_2, q_3$ on $\SU_{1, 1}$) is smooth as well. To prove that the function $t_{max}(p)$ is continuous at a point $p \in C, \ |p| = 0$, we need to verify that $\frac{\partial q^p_0}{\partial t}(t_{max}(p), p) \neq 0$
(by the implicit function theorem).

Assume (by contradiction) that $\frac{\partial q^p_0}{\partial t}(t_{max}(p), p) = 0$,
then by definition of $t_{max}(p)$ we have $q^p_0(t_{max}(p), p) = 0$. So, we get the following system of equations:
$$
\left\{
\begin{array}{l}
-\frac{\eta p_3}{2 I_1} (1+ \eta) \sin \frac{t \eta p_3}{2 I_1} - \frac{t \eta p_3^2}{4 I_1^2} \cos \frac{t \eta p_3}{2 I_1} = 0, \\
\cos \frac{t \eta p_3}{2 I_1} - \frac{t p_3}{2 I_1} \sin \frac{t \eta p_3}{2 I_1} = 0. \\
\end{array}
\right.
$$
Express $\cos \frac{t \eta p_3}{2 I_1}$ from the second equation and substitute it to the first one. From $\frac{\eta p_3}{2 I_1} \sin \frac{t \eta p_3}{2 I_1} \neq 0$ we get $(\frac{t p_3}{2 I_1})^2 = - \frac{1 + \eta}{\eta} < 0$, a contradiction.

2. Let us prove that $\lim_{\p \rightarrow 0}{\frac{2 \tau^h_0(|\p|) I_1}{|p|}} = +\infty$. Actually $\tau^h_0(\p)$ is the root of the equation $q^h_0(\p, \tau) = 0$. Note that $\sin (\tau\eta\p) \neq 0$, since otherwise $\cos (\tau\eta\p) = 0$. Dividing the equation by $\cosh \tau \sin (\tau\eta\p)$, we get the equation:
$$
\cot (\tau\eta\p) = \p \tanh \tau.
$$
Its first positive root lies in the interval $(-\frac{\pi}{2 \eta \p}, -\frac{\pi}{\eta \p})$. Thus, this root tends to infinity as $\p \rightarrow 0$. It is easy to see that
$$
|p| = \sqrt{\frac{I_1}{1 + \p^2\eta}} \rightarrow \sqrt{I_1},
$$
and the statement follows.
\qquad$\Box$
\medskip

We get the following description of the first Maxwell strata corresponding to the symmetry group $S$ in the pre-image and in the image of the exponential map.

{\Corollary
\label{crl-symmetries-first-maxwell-set}
$(1)$ When $\eta \leqslant -\frac{3}{2}$ we have
\begin{equation*}
\begin{array}{lll}
\M(S) & = & \{ (p, t) \in C^e \times \R_+ \ | \ t = \frac{2 \tau^e_0(|\p|) I_1}{|p|} \} \ \cup \\
      &   & \{ (p, t) \in C^p \times \R_+ \ | \ t = t^p_0(p) \} \ \cup \\
      &   & \{ (p, t) \in C^h \times \R_+ \ | \ \p \neq 0, \ \frac{2 \tau^h_0(|\p|) I_1}{|p|} \}, \\
\end{array}
\end{equation*}
$$
\Exp{(\M(S))} = Z := \Pi \{ q \in \SU_{1, 1} \ | \ q_0 = 0 \} \backsimeq \R^2
$$
is the plane of all central symmetries of the hyperbolic plane. \\
$(2)$ When $-\frac{3}{2} < \eta < -1$ we have
\begin{equation*}
\begin{array}{lll}
\M(S) & = & \{ (p, t) \in C^e \times \R_+ \ | \ |\p| > -\frac{3}{2\eta}, \ t = \frac{2 \tau^e_0(|\p|) I_1}{|p|} \} \ \cup \\
      &   & \{ (p, t) \in C^e \times \R_+ \ | \ |\p| \leqslant -\frac{3}{2\eta}, \ \p \neq \pm 1, \ t = \frac{2 \pi I_1}{|p|} \} \ \cup \\
      &   & \{ (p, t) \in C^p \times \R_+ \ | \ t = t^p_0(p) \} \ \cup \\
      &   & \{ (p, t) \in C^h \times \R_+ \ | \ \p \neq 0, \ \frac{2 \tau^h_0(|\p|) I_1}{|p|} \}, \\
\end{array}
\end{equation*}
$$
\Exp{(\M(S))} = Z \cup R_{\eta},
$$
where the interval
$$
R_{\eta} := \{ R_{0, \pm \varphi} \in \PSL_2(\R) \ | \ \varphi \in (-2\pi(1+\eta), \pi] \}
$$
consists of the rotations around the center of the Poincar\'{e} disk model. (The rotation around the center of the Poincar\'{e} disk model by the angle $\varphi$ is denoted by $R_{0, \varphi}$.)
}
\medskip

\emph{Proof.}
The statements about $\M(S)$ follow from Propositions~\ref{prop-comparsion-e-tau0-tau3}, \ref{prop-comparsion-e-tau0-pi}, \ref{prop-comparsion-h-tau0-tau3}, \ref{prop-comparsion-p-t0-t3}.
For description of $\Exp{(\M(S))}$ recall that geodesics reach the set $\Pi \{q \in \SU_{1, 1} \ | \ q_0 = 0\}$ at the time corresponding to the values $\tau^e_0(\p), t^p_0(p), \tau^h_0(\p)$. This set consists of central symmetries (see Section~\ref{section-model}). It remains to show that we can get any central symmetry in the image of the first Maxwell stratum. Use continuity of the exponential map and continuity of the first Maxwell time corresponding to symmetries of the exponential map (Lemma~\ref{lemma-firt-maxwell-time-continious}). Actually, for
$|\p| = 1$ when $\eta \leqslant -\frac{3}{2}$ or $|p_3| = -\frac{3}{2\eta}$ when $\eta > -\frac{3}{2}$
a corresponding geodesic at the first Maxwell time reaches a point with $q_1 = q_2 = 0$.
Because of $\tau^h_0(\p) \in (-\frac{\pi}{2 \eta \p}, -\frac{\pi}{\eta \p})$ we have
$\tau^h_0(\p) \rightarrow \infty$ as $p_3 \rightarrow 0$, it follows $(q_1(\tau^h_0(\p)))^2 + (q_2(\tau^h_0(\p)))^2 \rightarrow +\infty$. A continuous function $(q_1(\tau^h_0(\p)))^2 + (q_2(\tau^h_0(\p)))^2$ takes all values of the interval $[0, +\infty)$. We can achieve any direction of the vector $(q_1(\tau^h_0(\p)), q_2(\tau^h_0(\p)))$ by an appropriate choice of $p_1, p_2$.

For description of the stratum $R_{\eta}$ notice that the function $q_0(\pi) = -\cos{(\pi \eta \p)}$ is continuous at the interval $(1, -\frac{3}{2 \eta}]$. Thus, this function takes all values in the interval from its minimum to its maximum $[0, \cos{(\pi (1+\eta))})$. This corresponds to rotations around the center of the Poincar\'{e} disk model by angles $(-2 \pi(1 + \eta), \pi]$. Another interval (corresponding to the rotations in the opposite direction) is obtained by opposite values of $\p$.
\qquad$\Box$
\medskip

\section{\label{section-conjugate-time}Conjugate time}

{\Def
\emph{A conjugate point} is a critical value of the exponential map. \emph{A conjugate time}
is a time when a geodesic with arclength parametrization reaches the conjugate point.
}

{\Proposition
\label{prop-conj-time}
Consider a geodesic with an initial momentum $p \in C$.
For a time-like initial momentum and $\p \neq \pm 1$ there are two series of conjugate times:
$$
t_{2k-1} = \frac{2 I_1 \pi k}{|p|}, \qquad t_{2k} = \frac{2 I_1 \tau_k(p)}{|p|}, \qquad k \in \mathbb{N},
$$
where $\tau_k(p)$ is the $k$-th positive root of the equation
$$
\tan{\tau} = -\tau\eta\frac{1-\p^2}{1+\eta\p^2}.
$$
For $\p = \pm 1$ these two series merge to one series:
$$
t_k = \frac{2 I_1 \pi k}{|p|}, \qquad k \in \mathbb{N}.
$$
For light- or space-like initial momenta the corresponding geodesics have no conjugate points.
}
\medskip

\emph{Proof.}
Calculate the Jacobian of the exponential map. To make calculations more easy and independent of the type of an initial covector use notation~(\ref{def-trig-hyp}).

For time- and space-like covectors $p$ the Jacobian is equal to
$$
J(\tau, p) = s^2(\tau, p) \left(
\frac{\partial q_0}{\partial \p} \frac{\partial q_3}{\partial \tau} -
\frac{\partial q_3}{\partial \p} \frac{\partial q_0}{\partial \tau}
\right).
$$
The partial derivatives are equal to
\begin{equation}
\label{eq-diff}
\begin{array}{lcl}
\frac{\partial q_0}{\partial \p} & = &
-\tau \eta c(\tau, p) \se - s(\tau, p) \se - \tau \eta \p s(\tau, p) \ce,\\
\frac{\partial q_3}{\partial \p} & = &
\tau \eta c(\tau, p) \ce + s(\tau, p) \ce - \tau \eta \p s(\tau, p) \se,\\
\frac{\partial q_0}{\partial \tau} & = &
-\type(p) (1 + \type(p) \eta \p^2) s(\tau, p) \ce - \p (1 + \eta) c(\tau, p) \se,\\
\frac{\partial q_3}{\partial \tau} & = &
-\type(p)(1 + \type(p) \eta \p^2) s(\tau, p) \se + \p (1 + \eta) c(\tau, p) \ce.\\
\end{array}
\end{equation}
Substituting these expressions to the formula of Jacobian, we get
\begin{equation}
\label{eq-jcb-of-exp-map}
J(\tau, p) = \type(p) s^3(\tau, p) \left[
\tau \eta (1 - \type(p) \p^2) c(\tau, p) + (1 +  \type(p) \eta \p^2) s(\tau, p)
\right].
\end{equation}

Now we find positive roots of the function $J(\tau, p)$. For a time-like covector $p$ the first multiplier $s^3(\tau, p)$ equals $\sin^3{\tau}$ and vanishes at the points $\pi k, \ k \in \mathbb{N}$. For a space-like covector $p$ this multiplier is equal to $\sinh^3{\tau}$, and it has no positive roots.

Consider roots of the second multiplier of the Jacobian (the expression in the square brackets).

Note that $\eta < -1$ implies $\eta \p^2 < -\p^2$. Thus, $1 + \eta \p^2 < 1 - \p^2 \leqslant 0$ for
a time-like covector $p$ ($\p^2 \geqslant 1$). But $1 - \eta \p^2 > 1 + \p^2 \geqslant 1$ for all $\p \geqslant 0$. It follows that
$1 +  \type(p) \eta \p^2 \neq 0$. So, if $\tau$ is a zero of the second multiplier of the Jacobian, then $c(\tau, p) \neq 0$, since otherwise $s(\tau, p) = 0$ and these functions can not vanish simultaneously. Therefore, we need to investigate roots of the equation
$$
\frac{s(\tau, p)}{c(\tau, p)} = -\tau \eta \frac{1 - \type(p) \p^2}{1 +  \type(p) \eta \p^2}.
$$
It is easy to see that the coefficient of $\tau$ in the right-hand side is non-negative.

For a time-like covector $p$ this coefficient is less than $1$, since $\eta < -1$ and $-\eta + \eta \p^2 > 1 + \eta \p^2$. This means that the line with such slope does not intersect the branch of
the plot of the function $\tan{\tau}$ passing through the origin. For $\p \neq \pm 1$ we get
$\tau_k(p) \in (\pi k, \pi k + \frac{\pi}{2})$. For $\p = \pm 1$ we have $\tau_k(p) = \pi k, \ k \in \mathbb{N}$.

For a space-like covector $p$ we have $-\eta - \eta \p^2 > 1 - \eta \p^2$. This means that the slope of the line is greater than $1$. Thus this line intersects the plot of the function $\tanh{\tau}$ at the origin only. It follows that for a space-like initial covector the corresponding geodesics have no conjugate points.

Consider now geodesics with light-like initial momenta. We will show that there are no conjugate points. Apply the argument from~\cite{sachkov-elastics-conj}. By contradiction assume that for a light-like covector $p^0$ there is a finite conjugate time $t_{conj}(p^0)$. The conjugate points on the geodesic are isolated~(\cite{agrachev}). So, there exists $t_1 > t_{conj}(p^0)$ such that $t_1$ is not a conjugate time for $p^0$. Consider the continuous curve $p^s : [0, 1] \rightarrow C$ such that covectors $p^s$ are space-like for $s \in (0, 1]$ and $\lim_{s\rightarrow+0}{p^s} = p^0$. In paper~\cite{agrachev} it was shown that the number of conjugate points
(taking into account multiplicity) on the geodesic arc
$q^s(t) = \pi \circ e^{t \vec{H}} (p^s, \id), \ t \in [0, t_1]$, is equal to the Maslov index~\cite{arnold-index-maslova} of the path $l^s(t) = e^{-t \vec{H}}_* T^*_{q^s(t)} G$ in the Grassmanian of Lagrangian subspaces of $T_{(\id, 0)} T^*G$. Due to homotopic invariance of the Maslov index, the number of conjugate points on the geodesic arcs $q^0(t), \ t \in [0, t_1]$, and $q^1(t), \ t \in [0, t_1]$, are equal. There are no conjugate points on the arc $q^1(t), \ t \in [0, t_1]$, thus there are no conjugate points on the geodesic arc with the light-like initial momentum $p^0$ for $t \in [0, t_1]$. We get a contradiction.
\qquad$\Box$
\medskip

{\Def
\emph{The first conjugate time} $t_{conj}(p)$ is the time when the arclength parametrized geodesic with the initial momentum $p$ reaches the first conjugate point. We set it equal to infinity for geodesics that have no conjugate points.
}

{\Corollary
\label{crl-first-conj-time}
$$
t_{conj}(p) = \left\{
\begin{array}{lcl}
\frac{2\pi I_1}{|p|}, & \text{for} & \type{(p)} = 1,\\
+\infty, & \text{for} & \type{(p)} \leqslant 0.\\
\end{array}
\right.
$$
}
\medskip

\emph{Proof.} Follows from Proposition~\ref{prop-conj-time}.
\qquad$\Box$
\medskip

{\Remark
The function $t_{conj}: C \rightarrow (0, +\infty]$ is continuous. Actually, as a time-like
initial covector tends to a light-like one, the corresponding conjugate time
$t_{conj}(p) = \frac{2 \pi I_1}{|p|}$ tends to infinity, since $|p| \rightarrow 0$.
}
\medskip

\section{\label{section-cut-set}Cut locus}

Consider the Maxwell strata in the pre-image of the exponential map that correspond to the
symmetry group $S$. Denote the domain bounded by the closure of these Maxwell strata by
$$
U = \{(p, t) \in C \times \R_+ \ | \ 0 < t < t_{max}(p) \},
$$
where $t_{max}(p)$ is the first Maxwell time~(\ref{eq-maxwell-time}) for the symmetry group $S$.

Note that $U$ is an open subset of $C \times \R_+$, since it is the domain under the graph of a
continuous function (see Lemma~\ref{lemma-firt-maxwell-time-continious}).

{\Proposition
\label{prop-diffeomorphism}
The map $\Exp : U \rightarrow G \setminus (\Exp{(\overline{\M(S)})} \cup \{\id\})$ is a diffeomorphism.
}
\medskip

\emph{Proof.}
We use the Hadamard global diffeomorphism theorem~\cite{krantz-parks}: a proper non-degenerate smooth map of smooth connected and simply connected manifolds of same dimension is a diffeomorphism.

The manifolds $U$ and $G \setminus (\Exp{(\overline{\M(S)})} \cup \{\id\})$ both are three-dimensional and
connected, since both of them are homeomorphic to the punctured ball.
Indeed, the first one is the domain under
the graph of a continuous function on $C \backsimeq S^2$. The second one is the open punctured solid torus
$G \setminus \{\id\}$ without closure of the Maxwell set, that is the union of the open meridional disk $Z$ of
the torus and the interval $R_{\eta}$ (for $\eta \in (-\frac{3}{2}, -1)$). The result of the subtraction is homeomorphic to the punctured open ball. Consequently the both manifolds are simply connected.

The exponential map is non-degenerate on $U$ (there are no critical points in $U$). Indeed, $U$ is the domain under the graph of the first Maxwell time~(\ref{eq-maxwell-time}). The first Maxwell time for time-like initial covectors is $\frac{2 I_1}{|p|}\min{(\pi, \tau^e_0(|\p|))}$ and it is less than
or equal to the first conjugate time $\frac{2 \pi I_1}{|p|}$. For initial covectors of other types
the first conjugate time is infinite (Corollary~\ref{crl-first-conj-time}).

Now we prove that the map
$\Exp : U \rightarrow G \setminus (\Exp{(\overline{\M(S)})} \cup \{\id\})$ is proper, i.e., pre-image
of a compact set $K \Subset G \setminus (\Exp{(\overline{\M(S)})} \cup \{\id\})$ is compact (closed and bounded).

Assume that $\Exp^{-1}{(K)}$ is unbounded. Then there exists a sequence $(p^n, t^n) \in \Exp^{-1}{(K)}$ such that $t^n \rightarrow +\infty$. Since $p^n$ belongs to a compact set $C$,
there is a converging subsequence. So, we can assume $p^n \rightarrow p^0 \in C$. Clearly $\p^0 = 0$, since otherwise $t_{max}(p^0)$ is bounded and $t^n$ can not tend to infinity.
Because of $\p^0 = 0$ the covectors $p^n$ are space-like for numbers $n$ big enough
and their lengths $|p^n|$ are separated from zero.

For the images $\Exp{(p^n, t^n)} \in K$ in the coordinates $q^n_0, q^n_1, q^n_2, q^n_3$ we have
$$
q^n_0 \rightarrow +\infty, \qquad q^n_3 \rightarrow 0,
\qquad
(q^n_1)^2 + (q^n_2)^2  =  \sinh^2\left(\frac{t^n |p^n|}{2 I_1}\right)((p^n_3)^2 + 1) \rightarrow +\infty,
$$
but $K$ is bounded. We get a contradiction.

Assume now that $\Exp^{-1}{(K)}$ is not closed. Since it is bounded, there is a sequence
$(p^n, t^n) \in \Exp^{-1}{(K)}$ converging to $(p, t) \in \overline{U} \setminus \Exp^{-1}{(K)}$.

Then the sequence $\Exp{(p^n, t^n)} \in K$ converges to $\Exp{(p, t)} \in K$, since the map $\Exp$ is
continuous and the set $K$ is compact.

Hence, if $(p, t) \in U$, then $(p, t) \in \Exp^{-1}{(K)}$. We get a contradiction.

For $(p, t) \notin U$ (i.e., located at the boundary of $U$) we have $(p, t) \in \overline{\M(S)}$ or
$(p, t) \in C$ (i.e., $t=0$), since the sets $\overline{\M(S)}$ and $C$ are closed. Then $\Exp{(p, t)}$
belongs to $\Exp{(\overline{\M(S)})}$ or is equal to $\id$. We get a contradiction with
compactness of $K$.

So the set $\Exp^{-1}{(K)}$ is compact, thus the map $\Exp$ is proper. Thereby hypotheses of the Hadamard theorem are satisfied and the statement of the proposition is true.
\qquad$\Box$
\medskip

{\Theorem
\label{theorem-cut-time}
$(1)$ When $\eta \leqslant -\frac{3}{2}$ the cut time is
$$
t_{cut}(p) =
\left\{
\begin{array}{rcll}
\frac{2 I_1 \tau^e_0(\bar{p}_3)}{|p|}, & \text{for} & p \in C^e, & \\
t^p_0(p_3),                            & \text{for} & p \in C^p, & \\
\frac{2 I_1 \tau^h_0(\bar{p}_3)}{|p|}, & \text{for} & p \in C^h, & \p \neq 0, \\
+\infty, & \text{for} & \p = 0. & \\
\end{array}
\right.
$$
$(2)$ When $\eta > -\frac{3}{2}$ the cut time is
$$
t_{cut}(p) =
\left\{
\begin{array}{rcll}
\frac{2 I_1 \tau^e_0(\bar{p}_3)}{|p|}, & \text{for} & p \in C^e, & |\p| > -\frac{3}{2\eta}, \\
\frac{2 I_1 \pi}{|p|},                 & \text{for} & p \in C^e, & |\p| \leqslant -\frac{3}{2\eta}, \\
t^p_0(p_3),                            & \text{for} & p \in C^p, & \\
\frac{2 I_1 \tau^h_0(\bar{p}_3)}{|p|}, & \text{for} & p \in C^h, & \p \neq 0,\\
+\infty, & \text{for} & \p = 0. & \\
\end{array}
\right.
$$
}

{\Theorem
\label{theorem-cut-locus}
$(1)$ When $\eta \leqslant -\frac{3}{2}$ the cut locus is the plane $Z$ consisting of central symmetries. \\
$(2)$ When $\eta > -\frac{3}{2}$ the cut locus is a stratified manifold $Z \cup \overline{R_{\eta}}$,
where
$$
\overline{R_{\eta}} = \{ R_{0, \pm \varphi} \in \PSL_2(\R) \ | \ \varphi \in [-2\pi(1+\eta), \pi] \}
$$
is the interval consisting of some rotations around the center of the Poincar\'{e} disk model.
}
\medskip

\emph{Proofs of Theorems~$\ref{theorem-cut-time}$ and $\ref{theorem-cut-locus}$}
immediately follow from Proposition~\ref{prop-diffeomorphism}.
\qquad$\Box$
\medskip

\begin{figure}[h]
\caption{Cut locus in $\PSL_2(\R)$.}
\label{pic-cut-locus}
\medskip
     \begin{minipage}[h]{0.5\linewidth}
        \center{\includegraphics[width=0.45\linewidth]{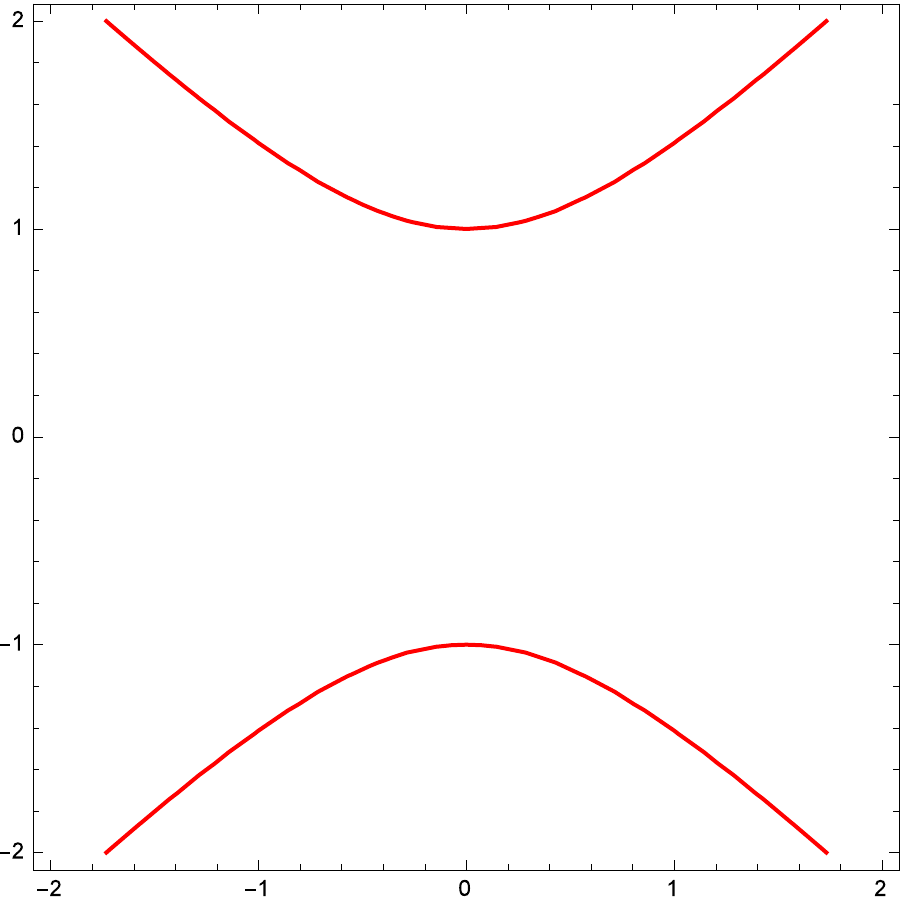} \\
        $\eta \leqslant -\frac{3}{2}, \ Z$}
     \end{minipage}
     \hfill
     \begin{minipage}[h]{0.5\linewidth}
        \center{\includegraphics[width=0.45\linewidth]{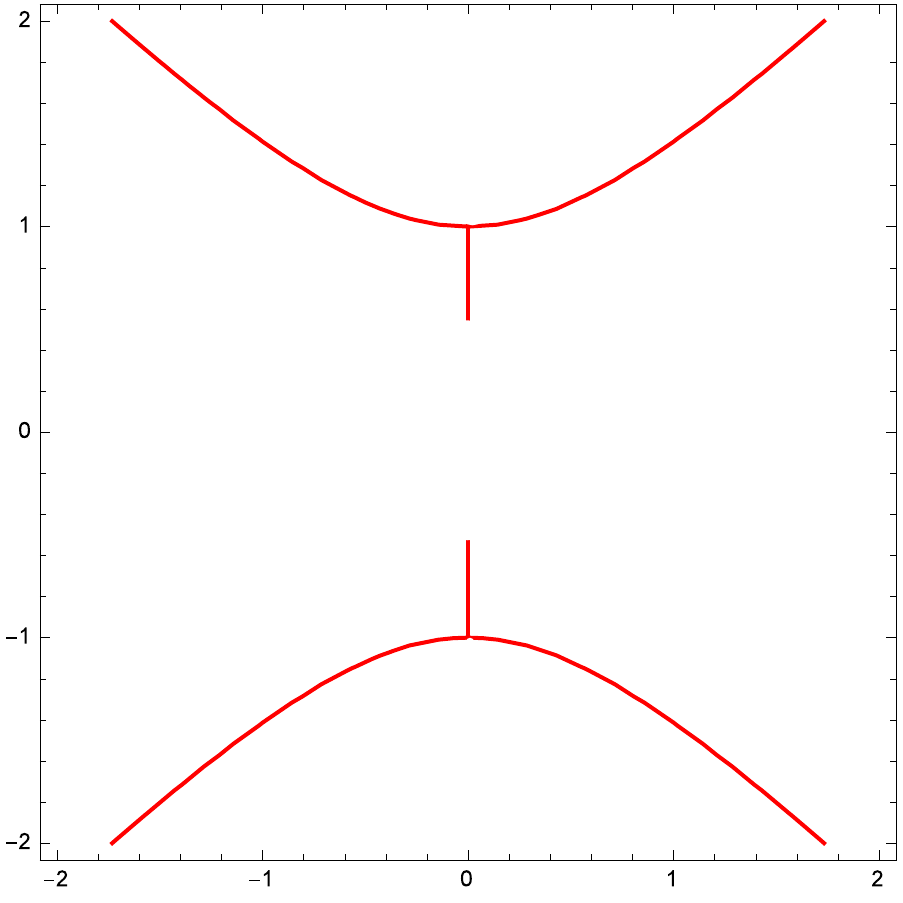} \\
        $\eta > -\frac{3}{2}, \ Z \cup \overline{R_{\eta}}$}
     \end{minipage}
\end{figure}

The cut locus is the surface of revolution of the contours presented in Figure~\ref{pic-cut-locus}
(in the model of $\PSL_2(\R)$ which is an open solid torus considered as the domain between two cups of a hyperboloid with the boundary identification).

\begin{figure}[h]
\caption{Wave front and cut locus for $\eta > -\frac{3}{2}$.}
\label{pic-wavefront}
\medskip
     \begin{minipage}[h]{0.18\linewidth}
        \center{\includegraphics[width=1\linewidth]{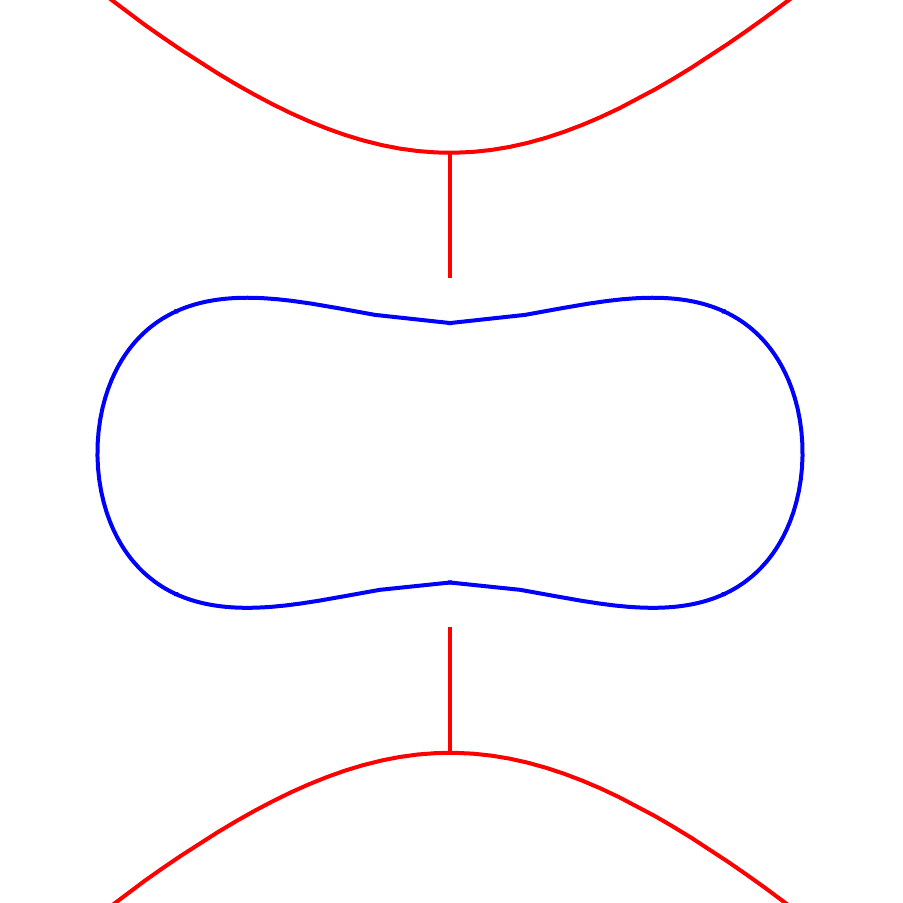}}
     \end{minipage}
     \hfill
     \begin{minipage}[h]{0.18\linewidth}
        \center{\includegraphics[width=1\linewidth]{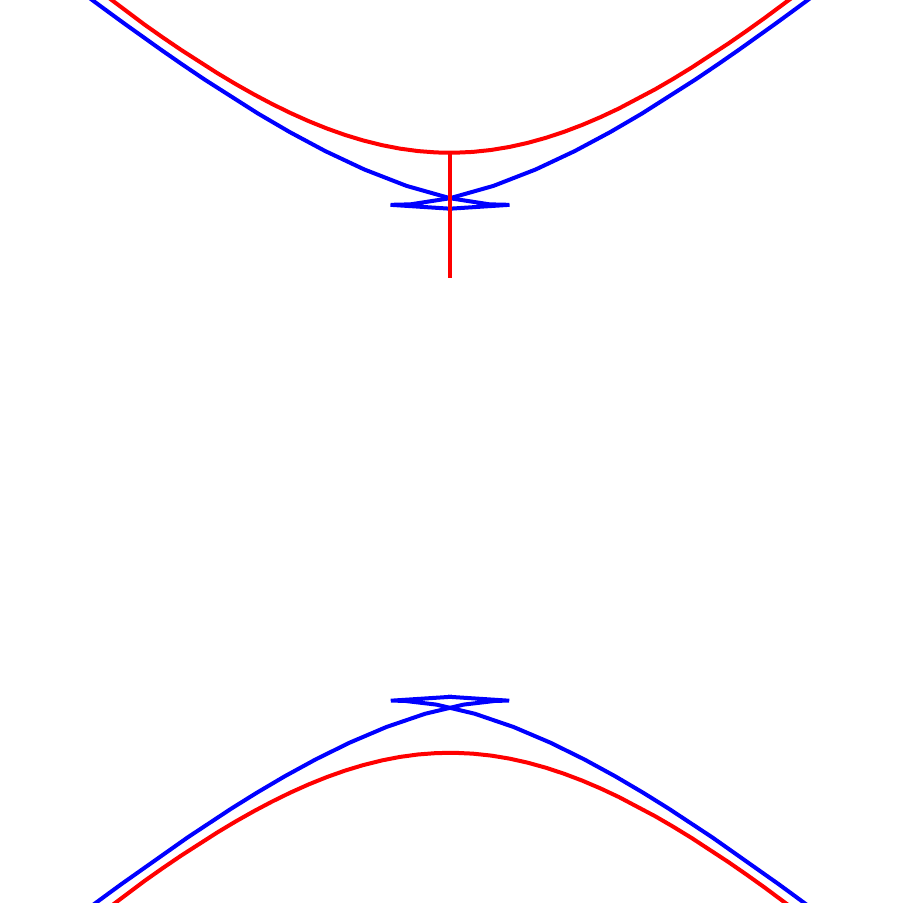}}
     \end{minipage}
     \hfill
     \begin{minipage}[h]{0.18\linewidth}
        \center{\includegraphics[width=1\linewidth]{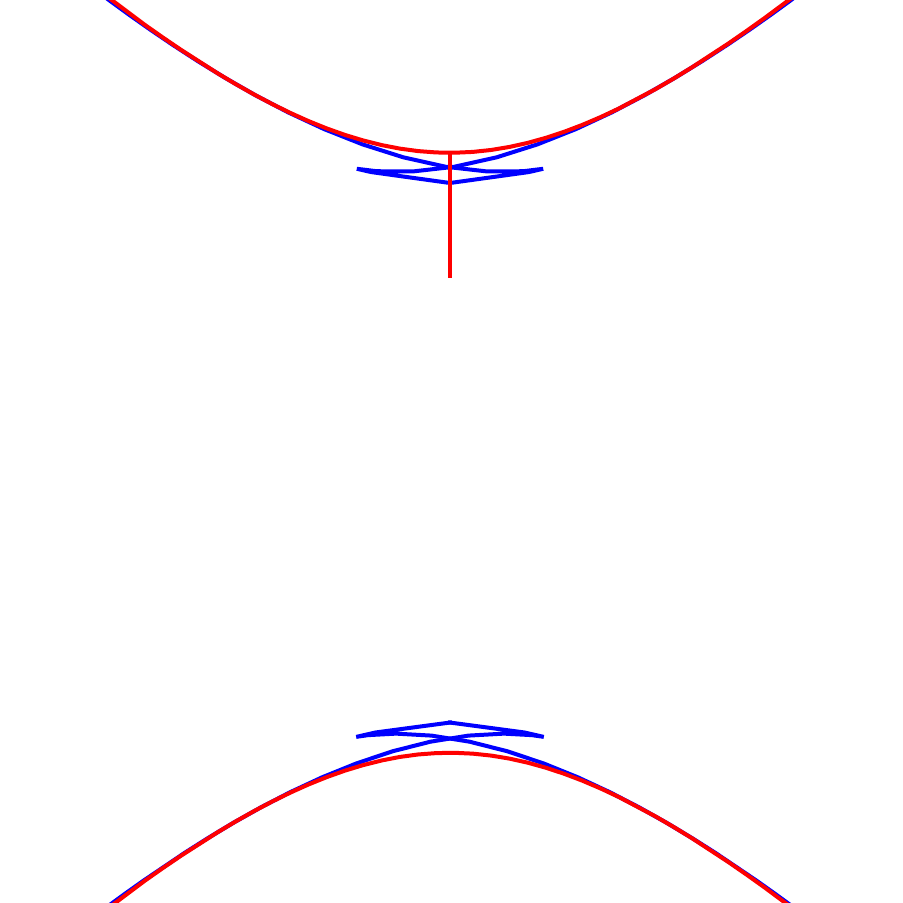}}
     \end{minipage}
     \hfill
     \begin{minipage}[h]{0.18\linewidth}
        \center{\includegraphics[width=1\linewidth]{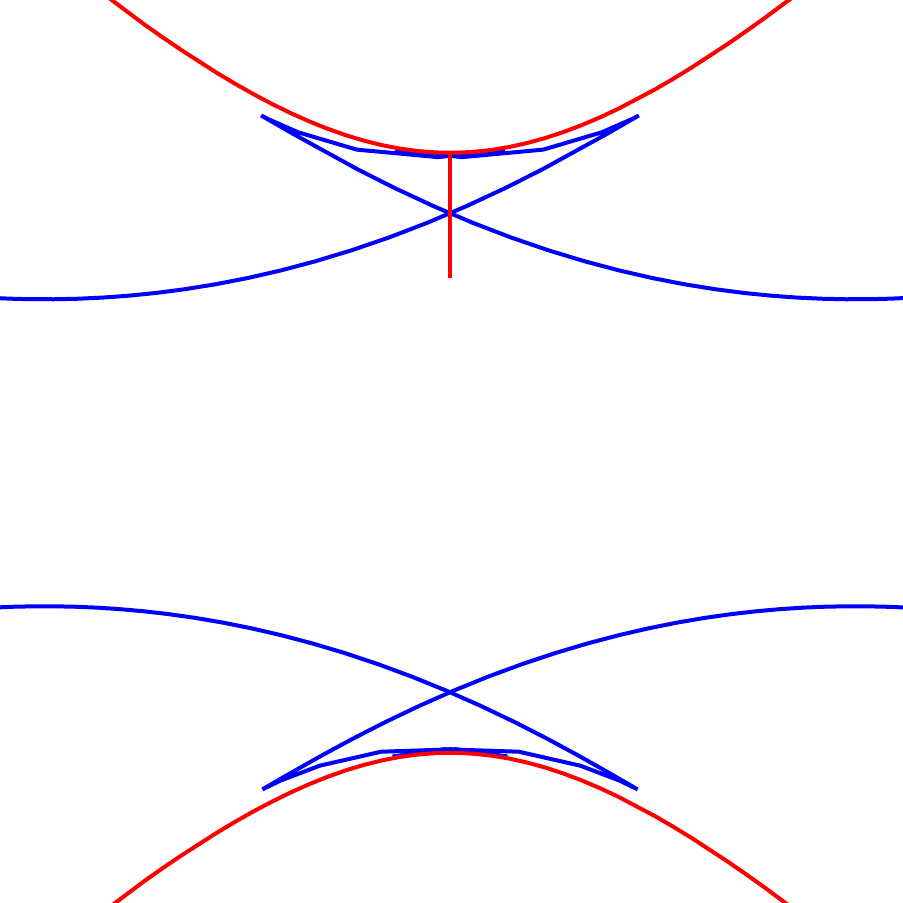}}
     \end{minipage}
     \hfill
     \begin{minipage}[h]{0.18\linewidth}
        \center{\includegraphics[width=1\linewidth]{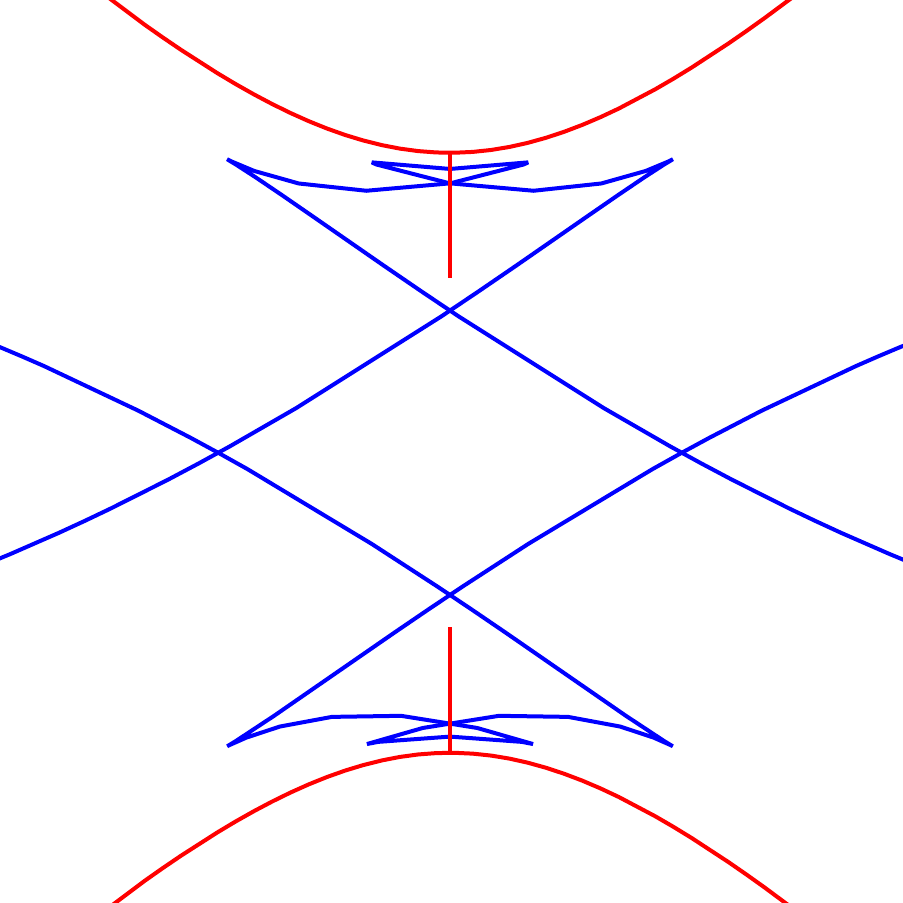}}
     \end{minipage}
\end{figure}

Propagation of the equidistant wave front is represented at Figure~\ref{pic-wavefront}
for $\eta > -\frac{3}{2}$.

{\Remark
\label{remark-maxwell-and-conjugate-points}
When $\eta \leqslant -\frac{3}{2}$ the cut locus coincides with the set of the first Maxwell points.
When $\eta > -\frac{3}{2}$ in the cut locus there are two conjugate points $R_{0, \pm 2 \pi (1 + \eta)}$ in
addition to the set of the first Maxwell points.
}
\medskip

\section{\label{section-injection-radius}Injectivity radius}

In this section we compute injectivity radius of the symmetric left-invariant Riemannian metric on $\PSL_2(\R)$. Recall that \emph{the injectivity radius} is the supremum of the set of numbers $T$ such that
the restriction of the exponential map to the set
$$
\{(p, t) \ | \ p \in C, \ 0 < t < T \}
$$
is injective.
It is clear that injectivity radius is equal to $\inf{\{t_{cut}(p) \ | \ p \in C\}}$.

Below we investigate the function $t_{cut}(p)$ defined on the sets $C^e$, $C^p$ and $C^h$, find and compare its local minima. The cut time is not a smooth function on $C$, but it is defined by the Maxwell times
corresponding to the strata $\M_0$ and $\M_{12}$, these times are smooth functions of the variable $\p$.

Denote the first Maxwell times corresponding to the strata $\M_0$ and $\M_{12}$ by
$$
t_0(\p) =  \frac{2I_1 \tau_0(\p)}{|p|}, \qquad
\tau_0(\p) = \left\{\begin{array}{lll}\tau^e_0(\p), &  p \in C^e, & \p \in [1, +\infty), \\
\tau^h_0(\p), &  p \in C^h, & \p \in (0, +\infty). \\
\end{array}\right.
$$
$$
t_{12}(\p)  = \begin{array}{lll} \frac{2\pi I_1}{|p|}, & p \in C^e, & \p \in [1, +\infty).\\ \end{array}
$$
Also introduce the functions:
$$
r(p, \eta) = 1 + \type{(p)}\eta\p^2,
$$
$$
j(\tau, p, \eta) = r(p, \eta)s(\tau, p)+\tau\eta(1-\type{(p)}\p^2)c(\tau, p).
$$

{\Proposition
\label{prop-diff-t0}
The next formulas are satisfied:
$$
|p| = \sqrt{\frac{I_1}{-\type{(p)}r(p, \eta)}}, \qquad \frac{\partial |p|}{\partial \p} = -\frac{\eta\p|p|}{\type{(p)}r(p, \eta)},
$$
$$
\frac{\partial t_0}{\partial \p} = -\frac{2I_1}{|p|} \frac{j(\tau, p, \eta)c(\tau, p)}{r(p, \eta)\p[\type{(p)}r(p, \eta)s^2(\tau, p) + (1+\eta)c^2(\tau, p)]}.
$$
}
\medskip

\emph{Proof.}
From $\frac{p_1^2}{I_1} + \frac{p_2^2}{I_1} + \frac{p_3^2}{I_3} = 1$ it follows $p_1^2 + p_2^2 = I_1 + (\eta + 1)p_3^2$. Thus,
$\Kil(p) = p_1^2 + p_2^2 - p_3^2 = I_1 + \eta p_3^2 = I_1 + \eta \p^2 |\Kil(p)| =
-\type{(p)}|\Kil{(p)}|$. Expressing $|\Kil{(p)}|$ and substituting it to $|p| = \sqrt{|\Kil{(p)}|}$, we get the first formula, the second one
can be produced just by computing derivative of the first one.

Next
$$
\frac{\partial t_0}{\partial \p} = 2I_1\frac{\frac{\partial \tau_0}{\partial \p}|p| - \tau_0 \frac{\partial |p|}{\partial \p}}{|p|^2}=
2I_1\frac{\frac{\partial \tau_0}{\partial \p}+\tau_0\frac{\eta\p}{\type{(p)}r(p, \eta)}}{|p|}.
$$
By the implicit function theorem we have
$$
\frac{\partial \tau_0}{\partial \p} = -\frac{\partial q_0}{\partial \p} \big/ \frac{\partial q_0}{\partial \tau}.
$$
Using expressions~(\ref{eq-diff}) of the partial derivatives of the function $q_0$, we get
$$
\frac{\partial \tau_0}{\partial \p} = -\frac{\tau \eta c(\tau, p) \se + s(\tau, p) \se + \tau \eta \p s(\tau, p) \ce}{\type(p) r(p, \eta) s(\tau, p) \ce + \p (1 + \eta) c(\tau, p) \se}.
$$
Consider the case $s(\tau, p) \se \neq 0$. Dividing numerator and denominator of the expression $\frac{\partial \tau_0}{\partial \p}$ by $s(\tau, p) \se$, we get
$$
\frac{\partial \tau_0}{\partial \p} = -\frac{\tau \eta \frac{c(\tau, p)}{s(\tau, p)} + 1 + \tau \eta \p \cot{(\tau\eta\p)}}{\type(p) r(p, \eta)\cot{(\tau\eta\p)}  + \p (1 + \eta) \frac{c(\tau, p)}{s(\tau, p)}}.
$$
Because of $q_0(\tau, p) = c(\tau, p) \ce - \p s(\tau, p) \se$, we obtain $\frac{c(\tau, p)}{s(\tau, p)} \cot{(\tau\eta\p)} = \p$. Next, $c(\tau, p) \neq 0$, since for $p \in C^h$ we have $c(\tau, p) = \cosh{\tau} \neq 0$, and for $p \in C^e$ if $c(\tau, p) = 0$, then $\p = 0$, in a contradiction
with $p \in C^e$.
Thus, $\cot{(\tau\eta\p)} = \p \frac{s(\tau, p)}{c(\tau, p)}$. Substituting it to the expression of
$\frac{\partial \tau_0}{\partial \p}$, we obtain
$$
\frac{\partial \tau_0}{\partial \p} = -\frac{\tau \eta c^2(\tau, p) + c(\tau, p)s(\tau, p) + \tau \eta \p^2 s^2(\tau, p)}{\p [\type(p) r(p, \eta) s^2(\tau, p)  + (1 + \eta) c^2(\tau, p)]}.
$$
Substituting this expression to the formula of $\frac{\partial t_0}{\partial \p}$, transforming to a common
denominator and using the equation $c^2(\tau, p) + \type{(p)} s^2(\tau, p) = 1$, we get the third formula of the proposition.

It remains to consider the case $s(\tau, p) \se = 0$.
Because of
$$
q_0(\tau, p) = c(\tau, p) \ce - \p s(\tau, p) \se = 0
$$
there are two cases.

The case $s(\tau, p) = 0$ and $\ce = 0$. Then
$$
\frac{\partial \tau_0}{\partial \p} = -\frac{\tau\eta}{\p(1+\eta)}, \qquad \frac{\partial t_0}{\partial \p} = -\frac{2I_1}{|p|} \frac{\tau\eta(1 - \type{(p)}\p^2)}{r(p, \eta) \p (1+\eta)},
$$
this coincides with the general formula.

The case $\se = 0$ and $c(\tau, p) = 0$. Then
$$
\frac{\partial \tau_0}{\partial \p} = -\frac{\tau\eta\p}{\type{(p)} r(p, \eta)}, \qquad \frac{\partial t_0}{\partial \p} = 0,
$$
this coincides with the general formula as well.
\qquad $\Box$
\medskip

{\Proposition
\label{prop-sgn-diff-t0}
The following equation is satisfied:
$$
\sgn{\left(\frac{\partial t_0}{\partial \p}\right)} = - \sgn{(\p)} \sgn{\left(j(\tau, p, \eta)\right)} \sgn{\left( c(\tau, p) \right)}.
$$
}
\medskip

\emph{Proof.}
In fact $r(\tau, \eta) = 1 + \type{(p)}\eta\p$. For $p \in C^e$ we have $\p \geqslant 1$. From $\eta < -1$ it follows $\eta \p^2 < -1$, i.e., $r(p, \eta) < 0$. For $p \in C^h$ we get $-\eta \p^2 > 0$ then
$r(p, \eta) > 0$. Hence, $\sgn{(r(p, \eta))} = -\type{(p)}$.

Therefore, the expression in the square brackets in denominator of the expression of $\frac{\partial t_0}{\partial \p}$ is negative. The statement of the proposition follows.
\qquad $\Box$
\medskip

{\Proposition
\label{prop-incr-decr-t0}
The function $t_0(\p)$ satisfies the properties: \\
$(1)$ it is increasing at the interval $[1, +\infty)$ for $p \in C^e$ when $\eta \leqslant -2$; \\
$(2)$ it is decreasing at the interval $[1, -\frac{2}{\eta}]$ and it is increasing at the interval
$[-\frac{2}{\eta}, +\infty)$ for $p \in C^e$ when $-2 < \eta \leqslant -\frac{3}{2}$; \\
$(3)$ it is decreasing at the interval $[-\frac{3}{2\eta}, -\frac{2}{\eta}]$ and it is increasing at the interval $[-\frac{2}{\eta}, +\infty)$ for $p \in C^e$ when $\eta > -\frac{3}{2}$; \\
$(4)$ it is decreasing at the interval $(0, +\infty)$ for $p \in C^h$.
}
\medskip

\emph{Proof.}
Notice that the expression $j(\tau, p, \eta)$ appears as a multiplier in
expression~(\ref{eq-jcb-of-exp-map}) of the Jacobian $J$ of the exponential map.
It was shown in the proof of Proposition~\ref{prop-conj-time} that for $p \in C^e$ the first positive zero of the function
$j(\tau, p, \eta)$ (of variable $\tau$) is greater than $\pi$, and for $p \in C^h$ this function has no
positive zeros.

(1--3) Computing $j(\pi, p, \eta)$, we have $-\pi\eta(1-\p^2) < 0$.
Thus, a continuous function $j(\tau, p, \eta)$ is negative for $\tau \in (0, \pi]$. From Proposition~\ref{prop-comparsion-e-tau0-pi} it follows that
the function $\tau_0^e(\p)$ is less than or equal to $\pi$ on the intervals $[1, +\infty)$ and $[-\frac{3}{2\eta}, +\infty)$ when $\eta \leqslant -\frac{3}{2}$ and $\eta > -\frac{3}{2}$ respectively. That is why $j(\tau, p, \eta)$ is negative under the hypotheses of the proposition.
This means that on the considered intervals the sign of $\frac{\partial t_0}{\partial \p}$ is equal to the sign of $\cos{\tau}$
due to Proposition~\ref{prop-sgn-diff-t0}.

It remains to determine when the sign of $\cos{\tau_0^e(\p)}$ changes, i.e.,
$\tau_0^e(\p) = \frac{\pi}{2}$. Let us prove that this happens at the point $\p = -\frac{2}{\eta}$
for $\eta > -2$.

Consider $\p > -\frac{2}{\eta}$, then $q_0^e(0, \p) = 1 > 0$ and on the other hand
$q_0^e(-\frac{\pi}{\eta\p}, \p) = -\cos{(-\frac{\pi}{\eta\p})} < 0$, since for $\p > -\frac{2}{\eta}$ the inequality $0 < -\frac{\pi}{\eta\p} < \frac{\pi}{2}$ holds. Hence, the function $q_0^e(\tau, \p)$ (of variable $\tau$, for $\p > -\frac{2}{\eta}$) at the endpoints of the interval
$[0, -\frac{\pi}{\eta\p}]$
have values of different signs. So, this continuous function has zero inside this interval.
Consequently, $\tau_0^e(\p) < \frac{\pi}{2}$ for $\p > -\frac{2}{\eta}$.

Consider now the case $\p < -\frac{2}{\eta}$. Our aim is to prove the inequality
$\tau_0^e(\p) > \frac{\pi}{2}$. Notice that this inequality is satisfied for $\p = 1$
(indeed, $\tau_0^e(1) = -\frac{\pi}{1+\eta} > \frac{\pi}{2}$). Assume that the inequality breaks at some
point of the interval $[1, -\frac{2}{\eta})$. Since the function $\tau_0^e(\p)$ is continuous, there exists $\hat{p}_3 \in [1, -\frac{2}{\eta})$ such that $q_0^e(\frac{\pi}{2}, \hat{p}_3) = -\hat{p}_3\sin{(\frac{\pi}{2}\eta\hat{p}_3)} = 0$. Solving this equation in the variable $\hat{p}_3$, we get $\hat{p}_3 \in \{\frac{2k}{\eta} \ | \ k \in \Z \}$. But this set does not intersect the interval
$[1, -\frac{2}{\eta})$, so we get a contradiction.

(4) The expression $j(\tau, p, \eta)$ is non-vanishing for $\tau > 0$ and $p \in C^h$. Substituting
$\tau = \frac{\pi}{2}$ to this expression, we get $r(p, \eta) > 0$ (see the proof of
Proposition~\ref{prop-sgn-diff-t0}). Thus, the sign of $\frac{\partial t_0}{\partial \p}$ is opposite to the sign of $\cosh{\tau}$, which is always positive. Consequently, the function $t_0(p)$
is decreasing for $p \in C^h$ and $\p \in (0, +\infty)$.
\qquad $\Box$
\medskip

{\Proposition
\label{prop-incr-decr-t12}
The function $t_{12}(\p) = \frac{2\pi I_1}{|p|}$ is increasing at the interval $[1, +\infty)$ for $p \in C^e$.
}
\medskip

\emph{Proof.}
Use the formula of $|p|$ from Proposition~\ref{prop-diff-t0}. The function $r(p, \eta) = 1 + \type{(p)}\eta\p^2$ is decreasing at the interval $[1, +\infty)$, because of $\eta < 0$. Thus, the function
$-\type{(p)} r(p, \eta)$ increases. Hence, $|p|$ decreases, the statement of the proposition follows.
\qquad $\Box$
\medskip

\begin{figure}[h]
\caption{Cut time for $p \in C^e$.}
\label{pic-injection-radius}
\medskip
     \begin{minipage}[h]{0.23\linewidth}
        \center{\includegraphics[width=1\linewidth]{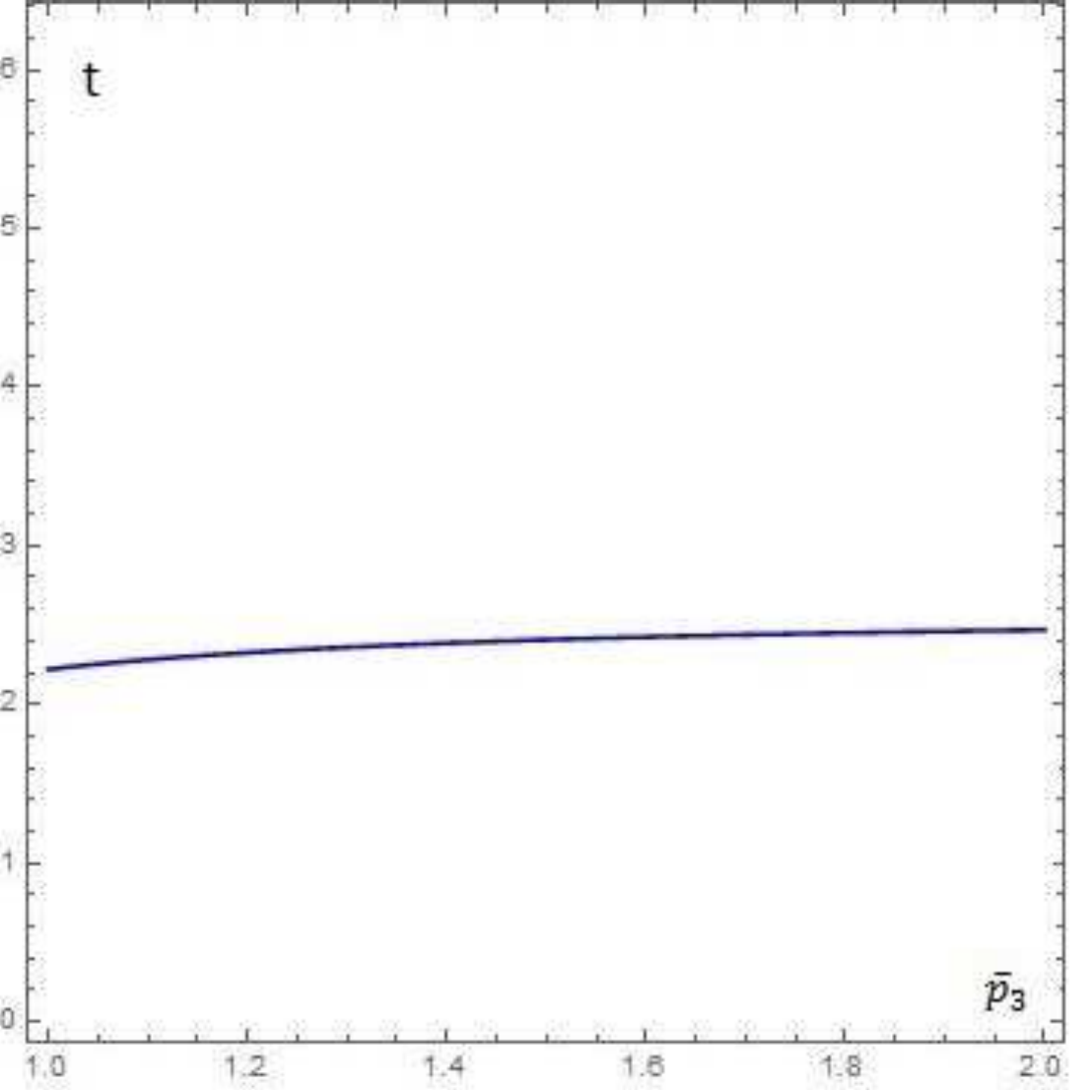}\\
        (a)\\ $\eta < -2$ }
     \end{minipage}
     \hfill
     \begin{minipage}[h]{0.23\linewidth}
        \center{\includegraphics[width=1\linewidth]{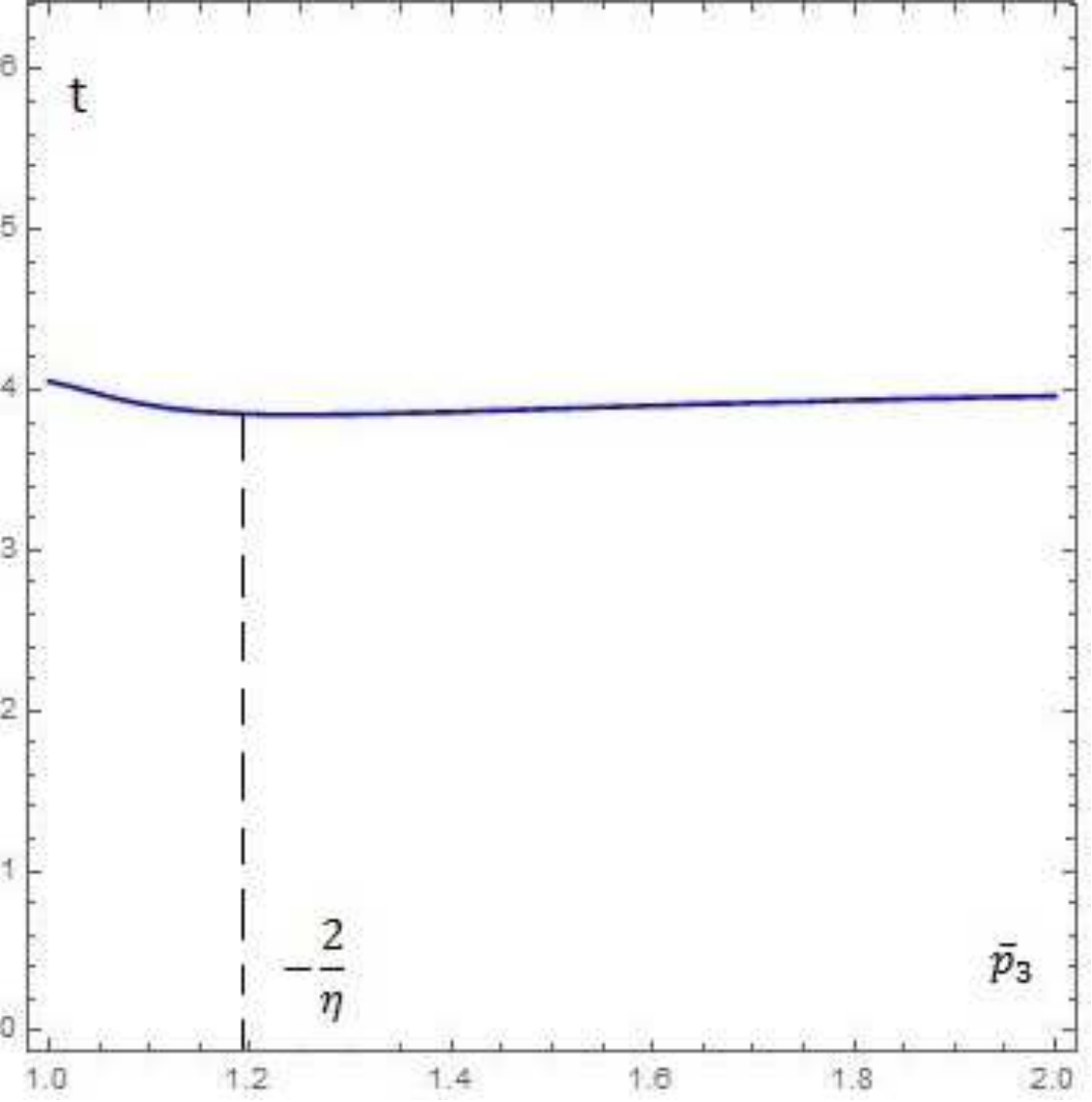}\\
        (b)\\ $-2 < \eta < -\frac{3}{2}$ }
     \end{minipage}
     \hfill
     \begin{minipage}[h]{0.23\linewidth}
        \center{\includegraphics[width=1\linewidth]{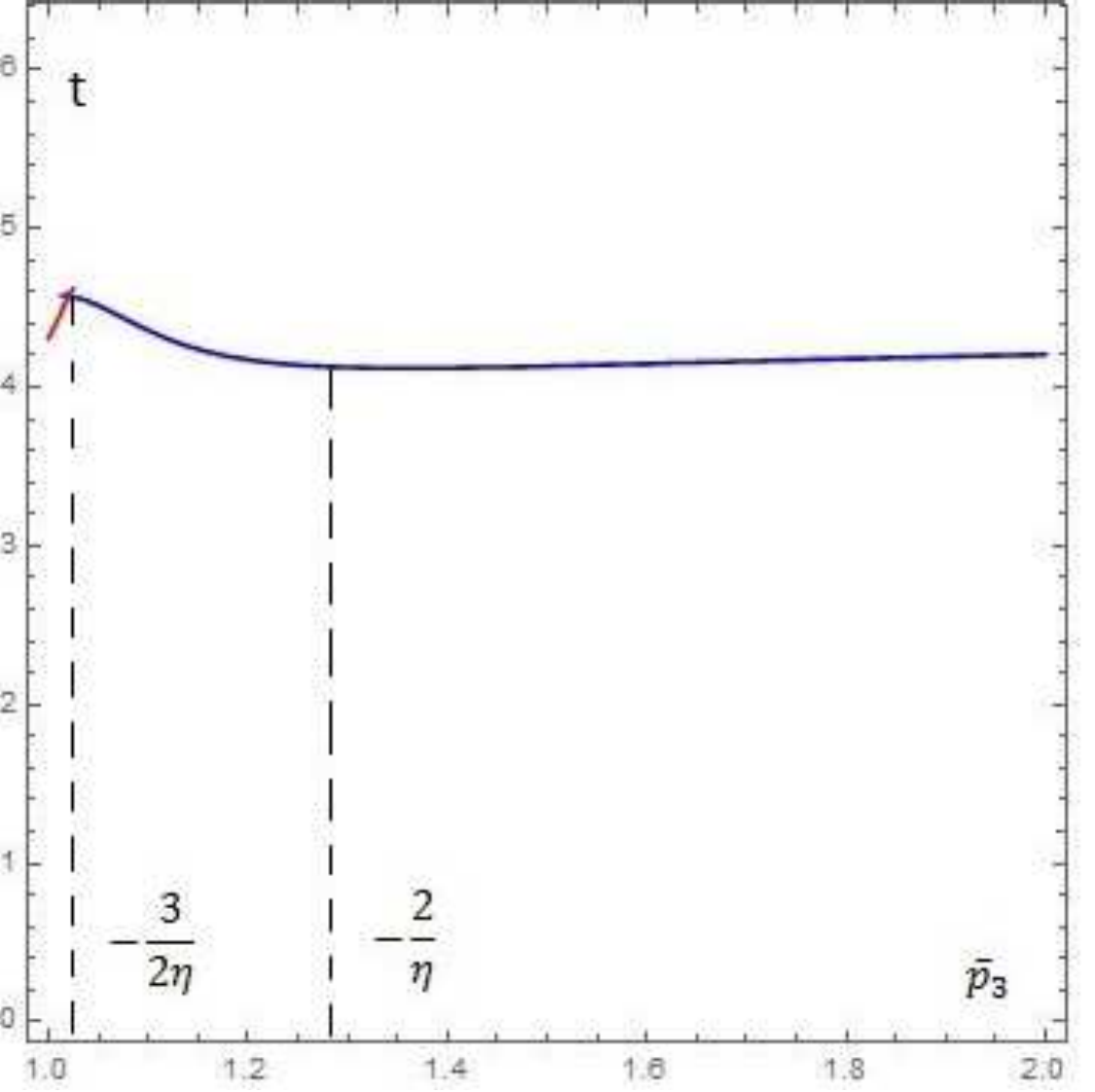}\\
        (c)\\ $-\frac{3}{2} < \eta < \frac{-3 - \sqrt{73}}{8}$ }
     \end{minipage}
     \hfill
     \begin{minipage}[h]{0.23\linewidth}
        \center{\includegraphics[width=1\linewidth]{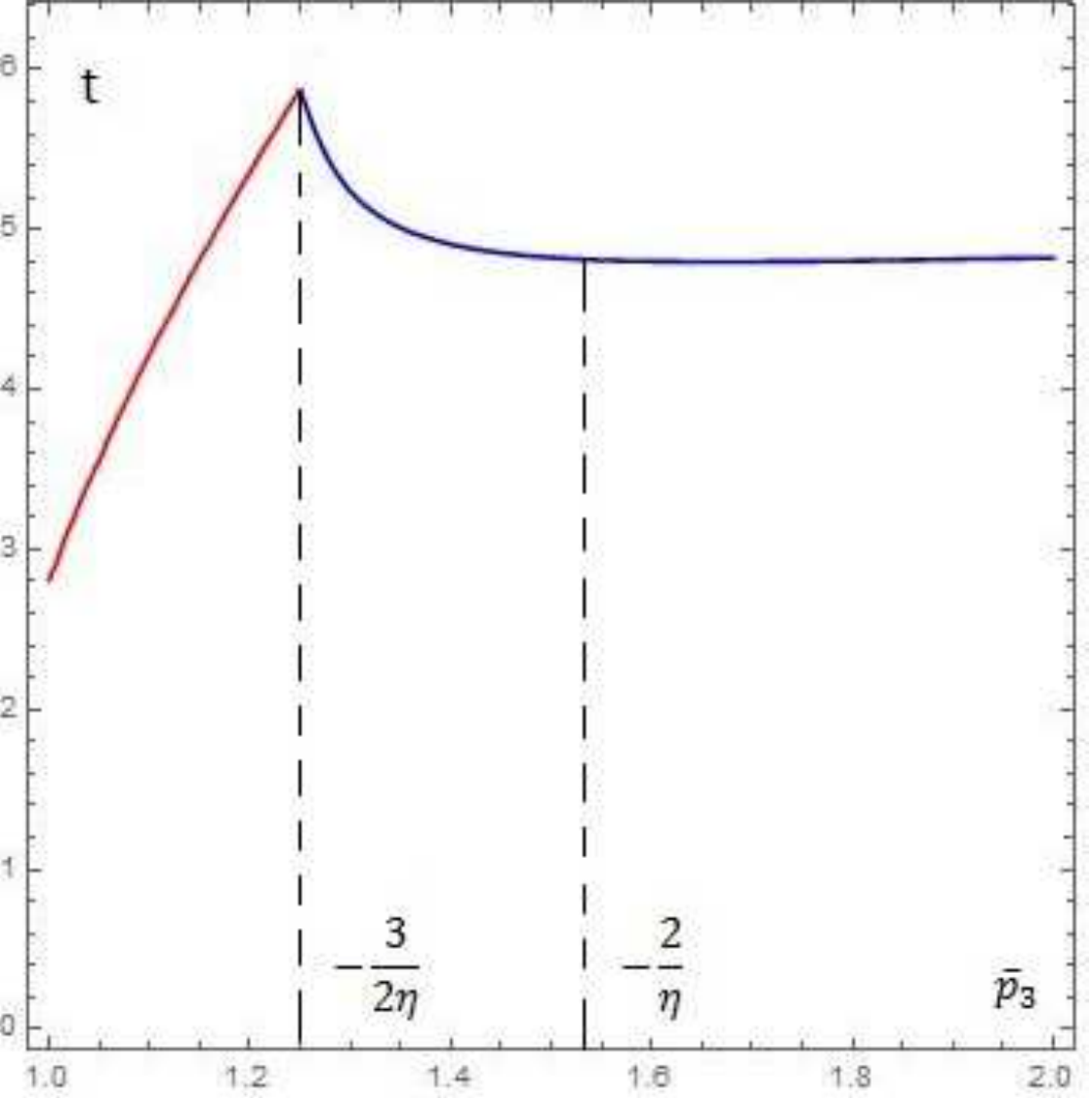}\\
        (d)\\ $\frac{-3 - \sqrt{73}}{8} < \eta < -1$ }
     \end{minipage}
\end{figure}

Figure~\ref{pic-injection-radius} presents plots of the cut time as the function of variable
$\p \in [1, +\infty)$ for time-like initial momenta and different values of the parameter $\eta$.

{\Corollary
\label{crl-inj-radius}
The injectivity radius of the symmetric left-invariant Riemannian metric on the group $\PSL_2(\R)$
is equal to \\
$(1)$ $\pi \sqrt{I_1} \sqrt{-\frac{1}{1+\eta}}$ when $\eta \leqslant -2$; \\
$(2)$ $\pi \sqrt{I_1} \sqrt{-\frac{\eta+4}{\eta}}$ when $-2 < \eta \leqslant \frac{-3-\sqrt{73}}{8}$; \\
$(3)$ $2\pi \sqrt{I_1}\sqrt{-(1+\eta)}$ when $\frac{-3-\sqrt{73}}{8} < \eta < -1$. }
\medskip

\emph{Proof.}
The injectivity radius is equal to the minimal value of the cut time $t_{cut}(p)$ for $p \in C$. From
Theorem~\ref{theorem-cut-time} it follows that the cut time is equal to the minimum of the Maxwell times corresponding to the strata $\M_0$ and $\M_{12}$.

From Propositions~\ref{prop-incr-decr-t0}, \ref{prop-incr-decr-t12} and continuity of the cut time
(Lemma~\ref{lemma-firt-maxwell-time-continious}) we obtain the following facts about the local minima of the
cut time.

($e$) There are two local minima of the cut time in the set $C^e$. They are the North $n$ and
the South $s$ poles of the ellipsoid $C$ (the points where $\p = \pm 1$).
The cut time has the same values at those points.
Besides, there are two circles of local minima $\{p \in C^e \ | \ \p = \pm \frac{2}{\eta} \}$.
The cut time is constant on these circles.
Denote by $m$ an arbitrary point of the circle in the North hemisphere, i.e., $\bar{m}_3 = -\frac{2}{\eta}$.

($p$) The cut time has no local minima in the set $C^p$ (for any point of $C^p$ there are an arbitrarily
close point of $C^e$ with a lower value of the cut time and an arbitrarily close point of $C^h$ with a greater value of the cut time).

($h$) In the set $C^h$ there are no local minima of the cut time (on the equator of the ellipsoid $C$
the cut time is infinite and it decreases along meridians from the equator to the poles).

The values of the cut time at the points of local minima are
$$
t_{cut}(n) = \left\{
\begin{array}{lll}
\frac{2\tau_0^e(\bar{n}_3)I_1}{|n|}, & \text{when} & \eta \leqslant -\frac{3}{2}, \\
\frac{2\pi I_1}{|n|}, & \text{when} & \eta > -\frac{3}{2}, \\
\end{array}
\right. \qquad
t_{cut}(m) = \frac{2\tau_0^e(\bar{m}_3)I_1}{|m|}.
$$
It is easy to see that for $\bar{n}_3 = 1$ the equation $q_0^e(\tau, \bar{n}_3) = \cos{\tau(1+\eta)} = 0$
has the first positive root $-\frac{\pi}{2(1+\eta)}$. Next
$|n| = \sqrt{\frac{I_1}{-(1+\eta)}}$. Thus
$$
t_{cut}(n) = \left\{
\begin{array}{lll}
\pi \sqrt{I_1} \sqrt{-\frac{1}{1+\eta}}, & \text{when} & \eta \leqslant -\frac{3}{2}, \\
2\pi \sqrt{I_1} \sqrt{-(1+\eta)}, & \text{when} & \eta > -\frac{3}{2}. \\
\end{array}
\right.
$$
Calculate now the value of the cut time at the point $m$. Note that $\tau_0^e(\bar{m}_3) = \frac{\pi}{2}$ and $|m| = \sqrt{-\frac{I_1\eta}{\eta+4}}$. Hence, $t_{cut}(m) = \pi \sqrt{I_1} \sqrt{-\frac{\eta+4}{\eta}}$.

Consider now different cases of the parameter $\eta$.

When $\eta \leqslant -2$ the cut time has no local minima (Proposition~\ref{prop-incr-decr-t0}) and
the injectivity radius is equal to $t_{cut}(n)$. We get case~(1), see Figure~\ref{pic-injection-radius}(a).

When $-2 < \eta \leqslant -\frac{3}{2}$ the cut time has a local minimum at the point $m$ which is the global minimum (Proposition~\ref{prop-incr-decr-t0}). We get case~(2), Figure~\ref{pic-injection-radius}(b).

When $\eta > -\frac{3}{2}$ we need to compare
$$
t_{cut}(n) = 2\pi \sqrt{I_1} \sqrt{-(1+\eta)} \qquad
\text{and} \qquad
t_{cut}(m) = \pi \sqrt{I_1} \sqrt{-\frac{\eta+4}{\eta}}.
$$
After elementary transformations
it is easy to see that $t_{cut}(n) < t_{cut}(m)$ if and only if $4\eta^2+3\eta-4 < 0$. This inequality is equivalent to $\eta \in (\frac{-3-\sqrt{73}}{8}, \frac{-3+\sqrt{73}}{8})$. It remains to use the inequalities:
$$
-\frac{3}{2} < \frac{-3-\sqrt{73}}{8} < -1 < \frac{-3+\sqrt{73}}{8}, \qquad \eta < -1.
$$
Thus when $\eta \in (-\frac{3}{2}, \frac{-3-\sqrt{73}}{8}]$ we have case~(2) presented in
Figure~\ref{pic-injection-radius}(c), and when $\eta \in (\frac{-3-\sqrt{73}}{8}, -1)$ we have case~(3), see Figure~\ref{pic-injection-radius}(d).
\qquad $\Box$
\medskip

{\Remark
The injectivity radius is a continuous function of the variable $\eta$.}
\medskip

\section{\label{section-sl2}Left-invariant Riemannian problem on $\SL_2(\R)$}

We use the same method of finding the cut locus as in the case of $\PSL_2(\R)$. Firstly, notice that the exponential map is described by formulas~(\ref{eq-time-geodesic}, \ref{eq-light-geodesic}, \ref{eq-space-geodesic}). Secondly, the symmetry group of the exponential map is the same that in the case of $\PSL_2(\R)$. The difference is that the set $\Exp{(\M_0)}$ is not a Maxwell stratum on $\SL_2(\R)$. In the case of $\PSL_2(\R)$ there are two geodesics that come to some point of this set at the same time,
but in the lift to $\SL_2(\R)$ these geodesics at that time are located in the different leaves of the covering
$\SL_2(\R) \rightarrow \PSL_2(\R)$.

The set of the first conjugate points and the first conjugate time are described  in the same way
as in the case of $\PSL_2(\R)$. Thus, for application of the Hadamard global diffeomorphism theorem
we need to compare the Maxwell time corresponding to the Maxwell strata $\Exp{(\M_{12})}$ and
$\Exp{(\M_3)}$ and the first conjugate time. The proposition below gives an answer for this question.

\begin{figure}[h]
\caption{Function $\tau^e_3$ and $\pi$.}
\label{pic-maxwell-time_sl2}
\medskip
     \begin{minipage}[h]{0.30\linewidth}
        \center{\includegraphics[width=1\linewidth]{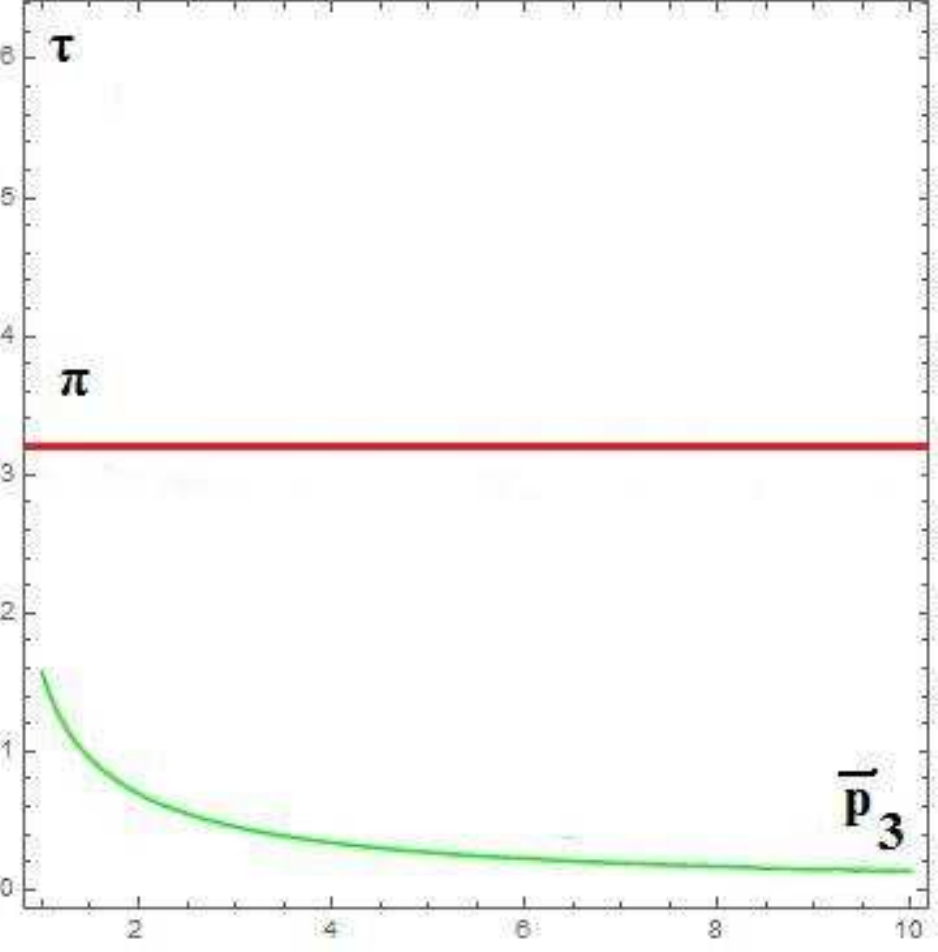} \\ $\eta < -\frac{3}{2}$}
     \end{minipage}
     \hfill
     \begin{minipage}[h]{0.30\linewidth}
        \center{\includegraphics[width=1\linewidth]{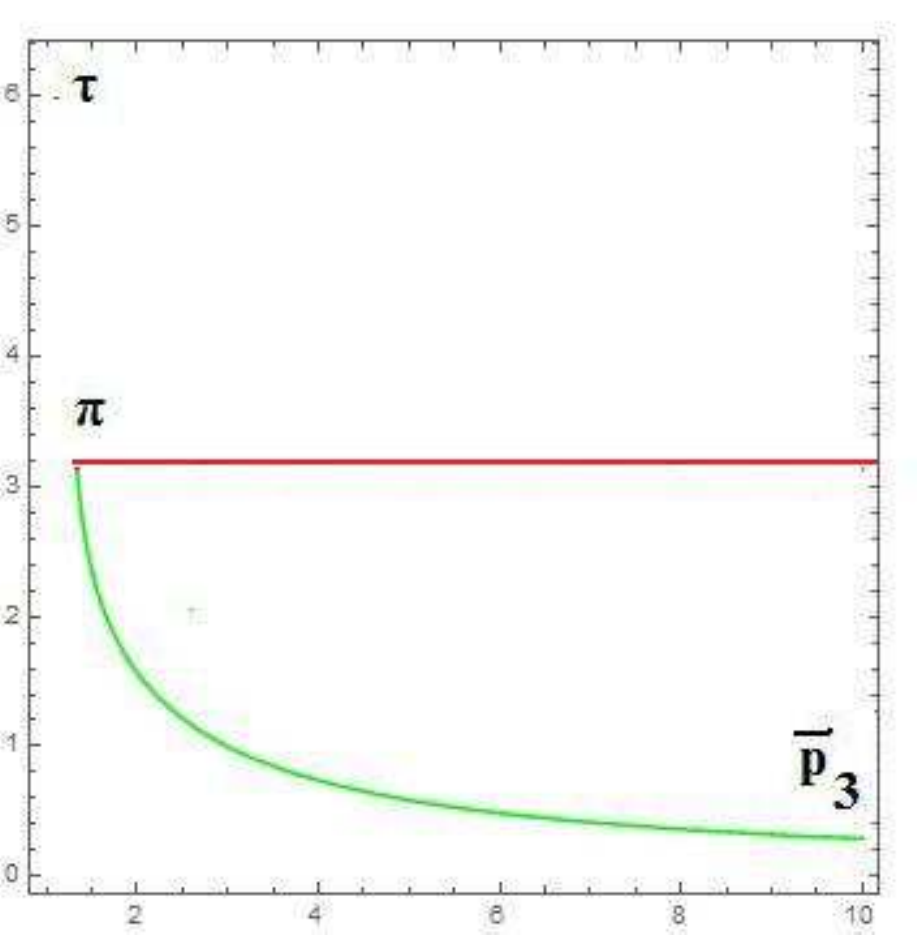} \\ $\eta = -\frac{3}{2}$}
     \end{minipage}
     \hfill
     \begin{minipage}[h]{0.30\linewidth}
        \center{\includegraphics[width=1\linewidth]{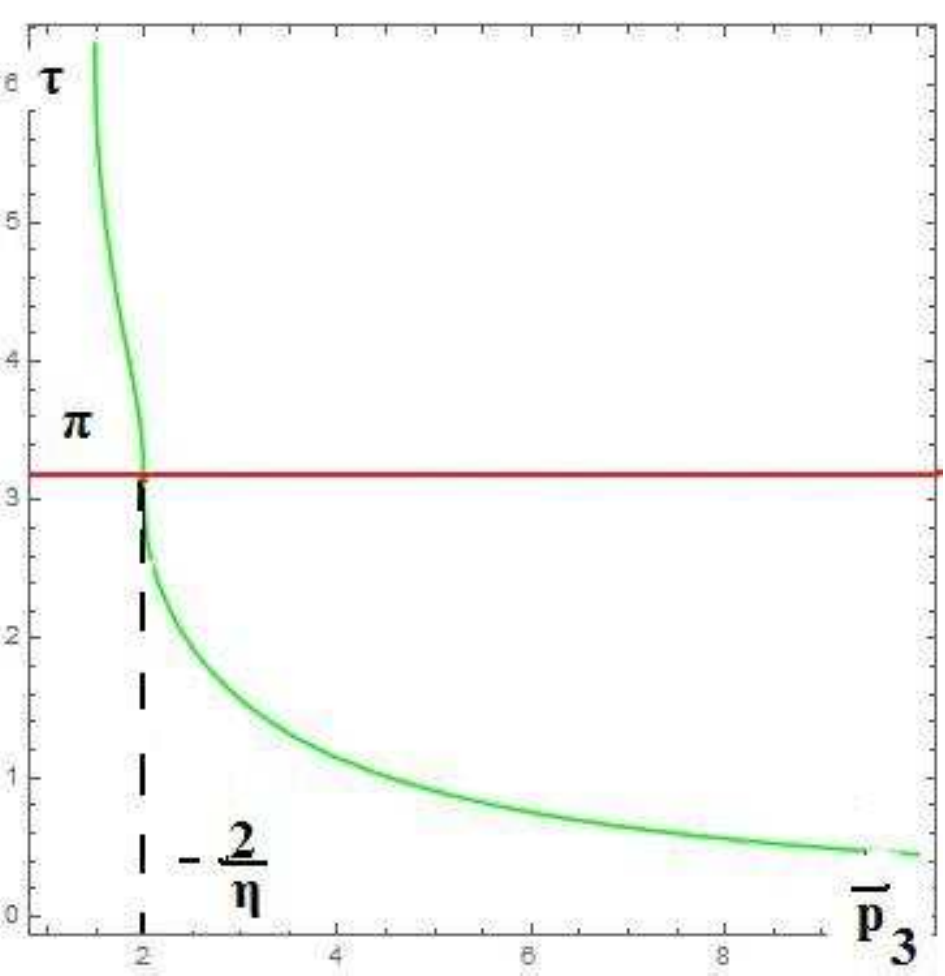} \\ $-\frac{3}{2} < \eta $}
     \end{minipage}
\end{figure}

{\Proposition
\label{prop-comparsion-tau3e_pi}
$(1)$ If $\eta \leqslant -\frac{3}{2}$, then for all $\p \in [1, +\infty]$ there holds the inequality $\tau^e_3(\p) \leqslant \pi$.\\
$(2)$ If $\eta > -\frac{3}{2}$, then
$\tau^e_3(\p) > \pi$ for $\p \in [1, -\frac{2}{\eta})$ and
$\tau^e_3(\p) \leqslant \pi$ for $\p \in [-\frac{2}{\eta}, +\infty)$.\\
See Figure~$\ref{pic-maxwell-time_sl2}$.
}
\medskip

\emph{Proof.}
(1) Note that $q_3^e(0) = 0$ and
$$
\frac{\partial q_3^e}{\partial \tau} = -(1+\eta \p^2) \sin{\tau} \sin(\tau \eta \p) + \p (1+\eta) \cos{\tau} \cos(\tau \eta \p)
$$
for $\tau = 0$ is equal to $\p (1 + \eta) < 0$.
The function $q_3^e$ of variable $\tau$ is differentiable. Thus, for all $\p \geqslant 1$ there exists
an arbitrarily small $\tau > 0$ such that $q_0^e(\tau, \p) < 0$. Hence, it is enough to find $\theta \in (0, \pi]$ such that $q_3^e(\theta) \geqslant 0$. Then, due to continuity of the function $q_3^e$ of variable $\tau$, there exists a root of $q_3^e$ at the interval $(0, \theta]$, i.e.,
a root that is less than or equal to $\pi$. Let us take
$$
\theta =
\left\{
\begin{array}{rl}
\frac{\pi}{2}, & \text{for} -\eta \p < 2, \\
-\frac{2 \pi}{\eta \p}, & \text{for} -\eta \p \geqslant 2. \\
\end{array}
\right.
$$
For $-\eta \p \geqslant 2$ we have $-\frac{2 \pi}{\eta \p} \leqslant \pi$. Next
$$
q_3^e(\theta) =
\left\{
\begin{array}{rl}
\p \cos (\frac{\pi}{2} \eta \p), & \text{for} -\eta \p < 2, \\
\p \sin(-\frac{2 \pi}{\eta \p}), & \text{for} -\eta \p \geqslant 2. \\
\end{array}
\right.
$$
In the first case $-\pi < \frac{\pi}{2} \eta \p \leqslant -\frac{3 \pi}{4}$, thus $\cos (\frac{\pi}{2} \eta \p) \geqslant 0$. In the second case $0 \leqslant -\frac{2 \pi}{\eta \p} \leqslant \pi$, it follows $\sin(-\frac{2 \pi}{\eta \p}) \geqslant 0$. Consequently $q_3^e(\theta) \geqslant 0$.

(2) First we prove the second part of the statement. For $\p \geqslant -\frac{2}{\eta}$
we have $-\eta \p \geqslant 2$. We can take the same $\theta$ as in the proof of statement~(1).

To prove the first part of statement~(2) notice that for $\p = 1$ we have $q_3^e(\tau) = \sin(\tau(1+\eta))$. The first positive root of this function is equal to
$-\frac{\pi}{1+\eta} > \pi$.
Assume (by contradiction) that there exists $\p' \in [1, -\frac{2}{\eta})$ such that $\tau_3^e(\p') < \pi$. The function $\tau_3^e$ is continuous, thus there exists $\widehat{p}_3 \in (0, \p')$ such that $q_3^e(\pi, \widehat{p}_3) = -\sin(\pi \eta \widehat{p}_3) = 0$. Hence $\widehat{p}_3 = \frac{k}{\eta}, \ k \in \Z$. It is clear that for all $k \in \Z$ the point $\widehat{p}_3$ lies outside of the interval $[1, -\frac{2}{\eta})$. We get a contradiction.
\qquad$\Box$
\medskip

For the symmetric left-invariant Riemannian problem on the group $\SL_2(\R)$ we describe below the cut
locus and the geometric interpretation of its image under the projection $\Pi$ onto the group of proper
isometries of the hyperbolic plane.

{\Theorem
\label{theorem-cut-locus-sl2}
$(1)$ When $\eta \leqslant -\frac{3}{2}$ the cut locus is the plane
$$
H := \{q \in \SU_{1, 1} \ | \ q_3 = 0 \},
$$
that maps to the plane of hyperbolic isometries corresponding to the sheafs of ultra-parallel lines
that are symmetric in the diameters of the Poincar\'{e} disk model.\\
$(2)$ When $\eta > -\frac{3}{2}$ the cut locus is a stratified manifold
$$
H \cup T_{\eta},
$$
where $T_{\eta} = \{ q = \pm (\cos(2\pi\p) + \sin(2\pi\p)k) \ | \ \p \in [1, -\frac{2}{\eta}]  \}$
is the interval that maps to the interval consisting of some rotations around the center of the Poincar\'{e}
disk model.
}
\medskip

\emph{Proof.} The proof is similar to the proof of Theorem~\ref{theorem-cut-locus}. This follows from Proposition~\ref{prop-comparsion-tau3e_pi} and the geometric interpretation of the subsets of the group of
proper isometries of the hyperbolic plane described in Section~\ref{section-model}.
\qquad$\Box$
\medskip

\section{\label{section-sub-riemannian}Connection with left-invariant sub-Riemannian\\ problem}

Identifying the Lie algebra $\mathfrak{g}$ with the space of pure imaginary split-quaternions, consider a decomposition
\begin{equation}
\label{eq-lie-algebra-decomposition}
\mathfrak{g} = \mathfrak{k} \oplus \mathfrak{p},
\end{equation}
where $\mathfrak{k} = \R k$ and $\mathfrak{p} = \R i \oplus \R j$.

Let $\Delta$ be the distribution of $2$-dimensional planes in $TG$ that is produced by the left shifts of the subspace $\mathfrak{p}$ of the Lie algebra. Endow the distribution $\Delta$ with the positive definite quadratic form $r_g(v) = (g^{-1}v, g^{-1}v)$, where
$g \in G$, $v \in \Delta_g = g \mathfrak{p}$ and $(\cdot, \cdot)$ is the Killing form.
Let $X_1, X_2$ be vector fields that form an orthonormal basis (with respect to the form $r_g$)
in the distribution $\Delta$ at every point.

Consider the following \emph{left-invariant sub-Riemannian problem}:
\begin{equation}
\label{eq-sub-riemannian-problem}
\dot{g} = u_1 X_1 + u_2 X_2, \qquad g(0) = \id, \qquad g(t_1) = g_1,
\qquad \frac{1}{2}\int_0^{t_1} (u_1^2 + u_2^2) \ dt \rightarrow \min.
\end{equation}

{\Theorem
\label{th-connection-with-sub-riemannian-problem}
For the left-invariant sub-Riemannian problem~$(\ref{eq-sub-riemannian-problem})$ on $\PSL_2(\R)$ $($or $\SL_2(\R))$ defined by decomposition~$(\ref{eq-lie-algebra-decomposition})$ and
the Killing form\\
$(1)$ the parametrization of geodesics, \\
$(2)$ the conjugate time, \\
$(3)$ the conjugate locus, \\
$(4)$ the cut time, \\
$(5)$ the cut locus \\
are produced from the same objects of the left-invariant Riemannian problem on $\PSL_2(\R)$
$($or $\SL_2(\R)$ respectively$)$ with $I_1 = I_2$ by passing to the limit
$I_3 \rightarrow \infty$.}

\begin{figure}[h]
\caption{Cut loci in sub-Riemannian and Riemannian cases.}
\label{pic-cut-locus-SR-R}
\medskip
     \begin{minipage}[h]{0.45\linewidth}
        \center{\includegraphics[width=0.5\linewidth]{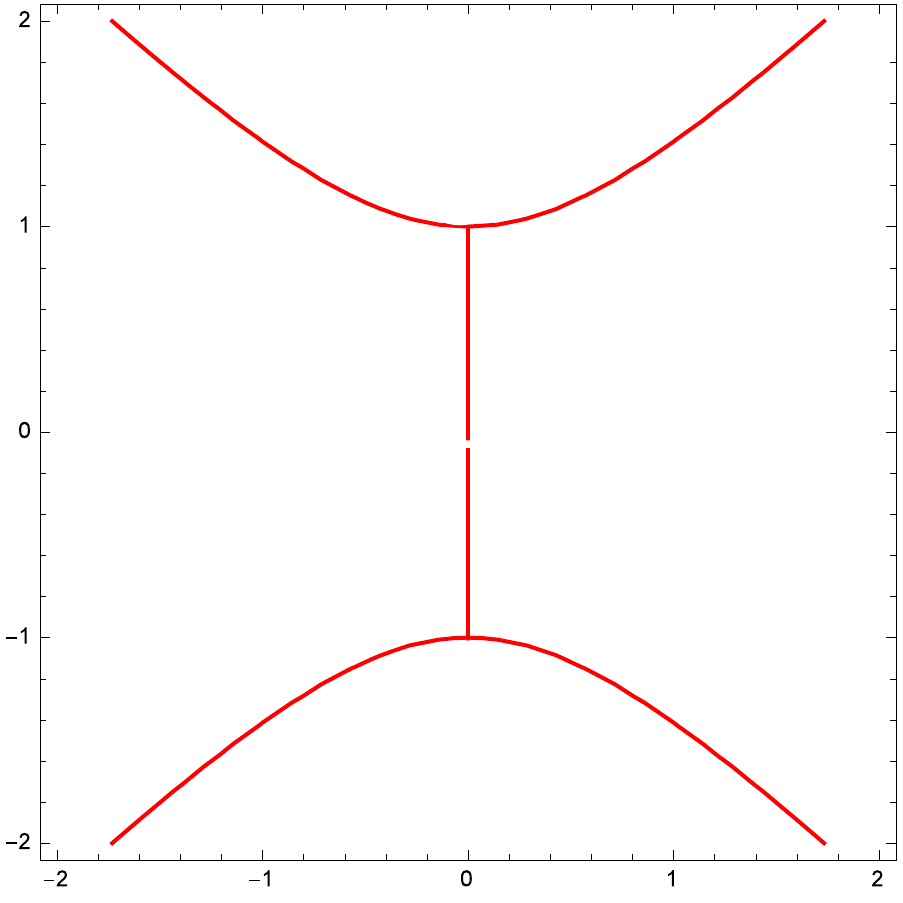} \\ sub-Riemannian metric}
     \end{minipage}
     \hfill
     \begin{minipage}[h]{0.45\linewidth}
        \center{\includegraphics[width=0.5\linewidth]{CutLocus_more_than_minus_3_2} \\ Riemannian metric, $\eta > -\frac{3}{2}$}
     \end{minipage}
\end{figure}

Figure~\ref{pic-cut-locus-SR-R} presents the cut loci for the sub-Riemannian and the Riemannian metrics
for $\eta > -\frac{3}{2}$ (the surfaces of revolution of the plotted contours).

\medskip
\emph{Proof.}
(1) Theorem~\ref{th-geodesics-parametrization} implies that the parametrization of geodesics on the considered
groups has the form
$$
g(t) = \exp \left(\frac{t}{I_1} p\right) \exp \left(\frac{t \eta p_3}{I_1} k\right),
$$
where $p \in \mathfrak{g}$ and $p = p_1 i + p_2 j + p_3 k$ is its split-quaternion representation.
As $\eta \rightarrow -1$ (this is equivalent to $I_3 \rightarrow \infty$) we get
$$
g(t) = \exp \left(\frac{t}{I_1}(p_1 i + p_2 j) + \frac{t}{I_1} p_3 k\right) \exp\left(-\frac{t}{I_1} p_3 k\right).
$$
This coincides with the known parametrization of sub-Riemannian geodesics (a proof could be found in
V.~Jurdjevic's book~\cite{jurdjevic}):
$$
g(t) = \exp (t(A_p + A_k)) \exp (-tA_k),
$$
where $A_k \in \mathfrak{k}$, $A_p \in \mathfrak{p}$, $r_{\id}(A_p) = 1$.
In V.~N.~Berestovskii's paper~\cite{berestovskii} a similar parametrization was got in Theorem~3:
$$
\gamma(t) = \exp{(t(\cos{\varphi_0}a + \sin{\varphi_0} b - \beta c))} \exp{(t \beta c)},
$$
where $\gamma$ is a geodesic, $a, b, c$ is a basis of the Lie algebra, the distribution is generated by
the vectors $a$ and $b$, the parameters $\varphi_0$ and $\beta$ define the initial covector.

For the sub-Riemannian problem on $\SL_2(\R)$ the same formula of parametrization of geodesics holds
(V.~N.~Berestovskii, I.~A.~Zubareva~\cite{berestovskii-zubareva-sl2}, Theorem~2).

The formulas relating the coordinates of initial covector in this paper and in
papers~\cite{berestovskii}, \cite{berestovskii-zubareva-sl2} are
$$
\begin{array}{cc}
|p|^2 = \beta^2 - 1, & \p = \frac{\beta}{\sqrt{|\beta^2 - 1|}}.\\
\end{array}
$$
\medskip

(2) The conjugate time for the sub-Riemannian problems on $\PSL_2(\R)$ and $\SL_2(\R)$ is finite
for $\beta > 1$ and is equal to $\frac{2 \pi}{\sqrt{\beta^2 -1}}$ (\cite{berestovskii}, \cite{berestovskii-zubareva-sl2}). In the Riemannian problem the conjugate time is finite only for
time-like initial covectors and it is equal to $\frac{2 \pi I_1}{|p|}$ (Proposition~\ref{prop-conj-time}), for $I_1 = 1$ it coincides with the sub-Riemannian conjugate time.

(3) The set of the first conjugate points is the circle $S^1 = \exp (\R k)$ both for the sub-Riemannian
and the Riemannian cases.

(4) The cut time in the sub-Riemannian problem on $\PSL_2(\R)$ was computed in~\cite{berestovskii} (Proposition~5). Below we give references (in parentheses) for the corresponding formulas from that paper.
For time-like initial covectors ($|\beta| > 1$) the cut time is equal to $\frac{2 \pi}{\sqrt{\beta^2 -1}}$
for $|\beta| > \frac{3}{\sqrt{5}}$ (52). For $1 < |\beta| \leqslant \frac{3}{\sqrt{5}}$ the cut time
is the first positive root of the equation (formulas~(54), (55)):
$$
-\cot{\frac{|\beta|t}{2}} = \frac{|\beta|}{\sqrt{\beta^2 - 1}} \tan{\frac{t\sqrt{\beta^2-1}}{2}}.
$$
For light-like initial covectors ($|\beta| = 1$) the cut time is the first positive root of the equation
(formulas~(50), (51)):
$$
-\cot{\frac{t}{2}} = \frac{t}{2}.
$$
For space-like initial covectors ($|\beta| < 1$) the cut time is the first positive root of the equation
(formulas~(48), (49)):
$$
-\cot{\frac{|\beta|t}{2}} = \frac{|\beta|}{\sqrt{1 - \beta^2}} \tanh{\frac{t\sqrt{1-\beta^2}}{2}}.
$$

Note that $\p = -\frac{3}{2 \eta}$ corresponds to $\beta = \frac{3}{\sqrt{5}}$. Thus,
the Riemannian cut time for light-like initial covectors for $\p > -\frac{3}{2 \eta}$ converges to the
cut time of the sub-Riemannian problem as $\eta \rightarrow -1$.

Clearly $\frac{|\beta|}{\sqrt{|\beta^2-1|}} = \p$, for $I_1 = 1$ we have
$\frac{|\beta|t}{2} = \tau \p$, $\frac{t\sqrt{\beta^2-1}}{2} = \tau$. Thus, for initial covectors of the other types
the equation $q_0 = 0$ converges to one of the equations above (depending on the type of initial covector).
Those equations and the equations $q^e_0 = 0, \ q^p_0 = 0, \ q^h_0 = 0$ for different values of $\eta$ do not have multiple roots. Hence, the first positive roots of the equations $q^e_0 = 0, \ q^p_0 = 0, \ q^h_0 = 0$
converge to the sub-Riemannian cut time as $\eta \rightarrow -1$.

For the sub-Riemannian problem on $\SL_2(\R)$, similar equations for the cut time were presented in
Theorem~6 of paper~\cite{berestovskii-zubareva-sl2}. Those equations are obtained from the
equations $q^e_3 = 0, \ q^p_3 = 0, \ q^h_3 = 0$ of the Riemannian cut time on $\SL_2(\R)$
(Theorem~\ref{theorem-cut-locus-sl2}) by passing to the limit $\eta \rightarrow -1$.
Note that for $|\beta| > \frac{2}{\sqrt{3}}$ the sub-Riemannian cut time is equal to
$\frac{2 \pi}{\sqrt{\beta^2 -1}}$. The initial covectors of such geodesics correspond to light-like $p$
with $\p > -\frac{2}{\eta}$.

(5) As $\eta \rightarrow -1$ the components $R_{\eta}$ and $T_{\eta}$ of the Riemannian cut loci on $\PSL_2(\R)$ and $\SL_2(\R)$ converge to the circle $S^1 = \exp (\R k)$ which is a component of the sub-Riemannian cut locus. The "global" part of the cut locus $Z$ ($H$ in case of $\SL_2(\R)$) is the
same for the Riemannian and the sub-Riemannian cases. The sub-Riemannian cut loci in $\PSL_2(\R)$ and
$\SL_2(\R)$ were described in papers of V.~N.~Berestovskii and I.~A.~Zubareva~\cite{berestovskii}, \cite{berestovskii-zubareva-sl2}, U.~Boscain and F.~Rossi~\cite{boscain-rossi}.
\qquad$\Box$
\medskip

\section*{\label{sec-appendix}Appendix. Some facts of hyperbolic geometry}

In this appendix we give some useful facts of the hyperbolic geometry.
Proofs can be found for example in book~\cite{prasolov-tikhomirov}.

{\Def
\emph{The Poincar\'{e} disk model} of the hyperbolic plane is the open unit disk
$\{z \in \mathbb{C} \ | \ |z| < 1 \}$. The boundary circle of the unit disk is called \emph{the absolute}. Points of the open unit disk are \emph{points} of the hyperbolic plane.
Consider Euclidean lines and circles that are orthogonal to the absolute. Arcs inside of the open unit disk are \emph{lines} of the hyperbolic plane. Clearly, there are infinite number of lines parallel
to a fixed line and passing through a point outside of that fixed line.
}

{\Def
\emph{The distance} $\rho(z_1, z_2)$ between two points $z_1$ and $z_2$ of the hyperbolic plane is defined as
$\rho(z_1, z_2) = \frac{c}{2} |\ln{|[u, v, z_1, z_2]|}|$, where
$u$ and $v$ are the intersection points of the line $z_1z_2$ and the absolute and
$$
[u, v, z_1, z_2] = \frac{z_1 - u}{z_1 - v} : \frac{z_2 - u}{z_2 - v}
$$
is \emph{the anharmonic ratio} of four points.
}

{\Remark
The parameter $c$ defines the eigenvalues $I_1 = I_2 = c$ of the Riemannian metric on the group of
proper isometries of the hyperbolic plane.
}

{\Theorem
$(1)$ Any proper isometry of the hyperbolic plane is determined by a M\"{o}bius transformation preserving the unit disk
$$
z \mapsto \frac{w_1 z + w_2}{\bar{w}_2 z + \bar{w}_1}, \qquad |w_1|^2 - |w_2|^2 = 1, \qquad
w_1, w_2 \in \mathbb{C}.
$$
$(2)$ Proper isometries form the group $\SU_{1, 1}$. \\
$(3)$ Any proper isometry is a composition of two reflections in lines. \\
$(4)$ There are three types of proper isometries: elliptic, parabolic and hyperbolic ones. The type is defined by the configuration of two lines. They can be intersecting, parallel one to another
(the intersection point belongs to the absolute) and ultra-parallel one to another (non-intersecting).\\
$(5)$ Orbits of these isometries are located on the curves that are orthogonal to the lines of elliptic, parabolic or hyperbolic sheaf respectively. Those curves are circle, oricircles or
equidistants respectively.
}

{\Remark
In the Poincar\'{e} half-plane model $\{z \in \mathbb{C} \ | \ \ImPart{z} > 0 \}$ of the hyperbolic plane
the group of proper isometries is the group of M\"{o}bius transformations of the form
$$
z \mapsto \frac{a z + b}{c z + d}, \qquad ad-bc = 1, \qquad a, b, c, d \in \R,
$$
that is isomorphic to $\PSL_2(\R)$. The transformation $z \mapsto i \frac{1 + z}{1 - z}$ maps the Poincar\'{e} disk model to the Poincar\'{e} half-plane model.
}

\end{document}